\documentclass{amsart}

\usepackage{amsmath}
\usepackage{amssymb}
\usepackage{amsthm}
\usepackage[english]{babel}
\usepackage{bm}
\usepackage[margin=0.5cm]{caption}
\usepackage{dsfont}
\usepackage{enumitem}
\usepackage[T1]{fontenc}
\usepackage{graphicx}
\usepackage[hidelinks]{hyperref} \usepackage[capitalise]{cleveref}
\usepackage{indentfirst}
\usepackage[utf8]{inputenc}
\usepackage{layout}
\usepackage{listings}
\usepackage{lscape}
\usepackage{mathrsfs}
\usepackage{mathtools}
\usepackage{nameref}
\usepackage[new]{old-arrows}
\usepackage{setspace}
\usepackage{soul}
\usepackage{stmaryrd}
\usepackage{subcaption}
\usepackage{textcomp}
\usepackage{thmtools}
\usepackage{xcolor}
\usepackage{xparse}

\usepackage{tikz}
\usepackage{tikz-cd}
\usetikzlibrary{shapes}

\title{On residual finiteness of graphs of free groups with cyclic edge groups}
\author{Adrien \textsc{Abgrall} and Zachary \textsc{Munro}}
\date{\today}

\DeclareMathOperator{\lk}{Link}

\begin{document}

\pagestyle{plain}
\frenchspacing
\parindent=15pt
\theoremstyle{plain}
\newtheorem{mthm}{Theorem}\renewcommand{\themthm}{\Alph{mthm}} 
\newtheorem{thm}{Theorem}[section]
\newtheorem{lm}[thm]{Lemma}
\newtheorem{pro}[thm]{Proposition}
\newtheorem{cor}[thm]{Corollary}
\newtheorem{conj}[thm]{Conjecture}
\theoremstyle{definition}
\newtheorem{defi}[thm]{Definition}
\newtheorem*{mdefi}{Definition}
\newtheorem{rem}[thm]{Remark}
\newtheorem*{ack}{Acknowledgments}
\theoremstyle{remark}
\newtheorem{q}[thm]{Question}
\newtheorem{ex}[thm]{Example}
\newtheorem{nota}[thm]{Notation}

\newcommand{\NN}{\mathbb{N}}
\newcommand{\ZZ}{\mathbb{Z}}
\newcommand{\QQ}{\mathbb{Q}}
\newcommand{\RR}{\mathbb{R}}
\newcommand{\CC}{\mathbb{C}}
\newcommand{\GG}{\mathcal{G}}
\renewcommand{\SS}{\mathbb{S}}
\newcommand{\TT}{\mathbb{T}}
\newcommand{\HH}{\mathbb{H}}
\newcommand{\FF}{\mathbb{F}}

\newcommand{\defeq}{\vcentcolon =}
\newcommand{\eqdef}{= \vcentcolon}
\newcommand{\dist}{\mathsf{d}}
\newcommand{\cat}{\mathrm{CAT}(0)}
\newcommand{\s}{\mathrm{Sep}}
\newcommand{\out}{\mathsf{Out}}
\newcommand{\aut}{\mathsf{Aut}}
\renewcommand{\emptyset}{\varnothing}
\renewcommand{\epsilon}{\varepsilon}

\newcommand{\n}[1]{\left|\left|\,#1\,\right|\right|}
\newcommand{\floor}[1]{\left\lfloor #1 \right\rfloor}
\newcommand{\ceil}[1]{\left\lceil #1 \right\rceil}
\newcommand{\gen}[1]{\left\langle #1 \right\rangle}
\newcommand{\card}[1]{\left | #1 \right |}

\newcommand{\bin}[1]{}
\newcommand{\todo}[1]{{\color{red}#1}}
\newcommand{\com}[1]{{\color{blue}#1}}

\begin{abstract}
    We characterize which groups splitting as finite graphs of free groups with cyclic edge groups are residually finite. Such a group $G$ is residually finite if and only if all its Baumslag-Solitar subgroups are residually finite. From a presentation of $G$, we construct a finite labeled graph $\Gamma$, and show that residual finiteness of $G$ is equivalent to an easily-detectable property of this graph. This characterization proves a conjecture of Wise from \cite{wise:separabilityGraphsFreeGroups}.
\end{abstract}

\maketitle

\section{Introduction}

A group $G$ is \emph{residually finite} if for every non-trivial $g\in G$ there exists a quotient to a finite group $G\to \overline{G}$ so that the image of $g$ in $\overline{G}$ is non-trivial. Residual finiteness implies an abundance of finite-index subgroups of $G$: an equivalent definition is that the intersection of all finite-index subgroups of $G$ is the trivial subgroup $1\in G$. The study of residual finiteness is classical in group theory, and many important classes of groups enjoy the property. Notably, free groups and surface groups are both residually finite \cite{hempel:residualFiniteness}.

Residual finiteness is also related to algorithmic problems in group theory. By an observation of McKinsey (and independently, Mostowski), finitely presented, residually finite groups have solvable word problem \cite{mckinsey:decisionProblem, mostowski:decidability}. Finitely presented, residually finite groups are also \emph{Hopfian}, i.e., every surjective endomorphism is injective. Baumslag and Solitar constructed the first examples of finitely generated, one-relator groups which are non-Hopfian \cite{baumslagSolitar:theOGpaper}. Their family of examples are now referred to as \emph{Baumslag-Solitar groups}. Studying the residual properties of Baumslag-Solitar groups and their generalizations has been a popular theme in group theory over the past few decades. For related discussions, see \cite{meskin:nonresiduallyFinite, moldavanskii:residualProperties}. 

In this paper, we study the family of groups which split as finite graphs of free groups with cyclic edge groups. This family of groups contains many interesting examples --- more than may initially be apparent. For example, free groups, orientable surface groups, and Baumslag-Solitar groups all split as finite graphs of free groups with cyclic edge groups. This family of groups has served as an important testing ground for conjectures and themes of exploration in geometric group theory. For example, Hsu and Wise characterized which of these groups are cubulated \cite{hsuWise:cubulating}, Wilton proved the surface subgroup conjecture in this family \cite{wilton:essentialSurfaces}, and Shepherd-Woodhouse proved quasi-isometric rigidity of a large sub-family \cite{shepherdWoodhouse:QIrigidity}.

Our main theorem is the following.
\begin{mthm}
\label{mthm:A}
    Suppose $G$ splits as a finite graph of free groups with cyclic edge groups. Then $G$ is residually finite if and only if the Baumslag-Solitar subgroups of $G$ are residually-finite.
\end{mthm}

Previously, Wise characterized the \emph{LERF} (i.e.~ subgroup separable) groups in this family, a much stronger property than residual finiteness \cite{wise:separabilityGraphsFreeGroups}. He proved that a group $G$ splitting as a finite graph of free groups with cyclic edge groups is LERF if and only if there exists no \emph{unbalanced} element, an infinite-order $g$ such that $g^n$ and $g^m$ are conjugate for $|n|$, $|m|$ distinct. We prove (Proposition~\ref{vuimpliesbs}) that this is equivalent to all the Baumslag-Solitar subgroups of $G$ being LERF.

Our proof uses an analogous characterization of residual finiteness: a group splitting as a graph of free groups with cyclic edge groups is residually finite if and only if it contains no \emph{very unbalanced elements} (see Definition~\ref{def:veryUnbalanced}). The existence of a very unbalanced element is easily detectable from the graph of groups splitting of $G$ (Remark~\ref{rem:algorithm}). In particular, there exists an algorithm which either confirms $G$ is residually finite or outputs a very unbalanced element. 

Theorem~\ref{mthm:A} was shown to hold for Generalized Baumslag-Solitar groups (i.e. finite graphs of infinite cyclic groups with infinite cyclic edge groups) by Levitt \cite[Corollary~7.7]{levitt:quotients}. The only cases where such a group is residually finite are when it is isomorphic to $BS(1,n)$ or does not contain unbalanced elements at all, in which case it is virtually $F\times \ZZ$ with $F$ free.

We use a topological approach to residual finiteness, realizing our graph of groups as the fundamental group of a graph of spaces $X$. Then $\pi_1 X$ is residually finite if and only if for each non-trivial $p\in \pi_1 X$ there exists a finite-degree cover $\widehat X\to X$ so that a closed path representing $p$ in $X$ lifts to a non-closed path in $\widehat X$. For a topological proof that the Baumslag-Solitar groups $BS(1,q)$ are residually finite, see Example~\ref{ex:BS(1,q)}. In the proof of Theorem~\ref{mthm:A}, we construct a sufficiently rich collection of finite-degree covers to witness residual finiteness. One of our tools is a variant of the Malnormal Special Quotient Theorem (MSQT) of Wise \cite{wise:structureQuasiconvexHierarchy}, which allows us to take controlled covers of the vertex spaces of $X$. Then we use small-cancellation theory to extend these covers to a covering of the whole $X$ and show that particular paths in $X$ must lift to non-closed paths there. The paper is organized as follows.

In Section~2 we summarize relevant background material. We first cover combinatorial complexes, paths and cycles, and a link criterion to verify that a map is a covering. We introduce graphs of spaces and their associated graphs of groups. We give several equivalent definitions of residual finiteness. We end with a brief review of small-cancellation theory, including van Kampen's Lemma, Greendlinger's Lemma, and a corollary of work by Agol and Wise. 

In Section~3 we define and discuss very unbalanced group elements. In general, a residually finite group can contain very unbalanced elements (Remark~\ref{rem:rfwithvu}). However, a residually finite group splitting as a graph of free groups cannot contain them, as they imply the existence of non residually finite Baumslag-Solitar subgroups (Proposition~\ref{vuimpliesbs}). In other words, when discussing graphs of free groups, very unbalanced elements serve as a certificate of non residual finiteness. This provides one direction of Theorem~\ref{mthm:A}. 

In Section~4 we define a labeled graph $\Gamma$ associated to a graph of free groups with cyclic edge groups $G$. We show that the existence (and non-existence) of a very unbalanced element in $G$ can be easily determined by structure of $\Gamma$. In particular, since $\Gamma$ is easily constructed from a standard presentation of $G$, there exists an algorithm determining the existence or non-existence of a very unbalanced element in $G$.

In Section~5 we prove that a group splitting as a finite graph of free groups with cyclic edge groups is residually finite if it contains no very unbalanced elements (Proposition~\ref{mainprop}), the other direction of Theorem~\ref{mthm:A}. Our proof is topological. We begin by defining \emph{normalized} representatives for elements in our group, as paths in the graph of spaces which will eventually lift to non-closed paths in the appopriate cover. Then we construct such finite-degree covers, which will witness residual finiteness. The structure of the graph $\Gamma$, defined in Section~4, guides the construction of the covers. We end with a case analysis of possible normalized representatives, showing that there is always a non-closed lift to one of our finite-degree covers.

\begin{ack}
    The authors thank the Centre de Recherches Mathématiques de Montréal (CRM) for hosting the \textit{Geometric group theory} thematic semester in 2023 where part of this work was completed. Adrien Abgrall thanks ANR-22-CE40-0004 GoFR for support.
    
    We are very thankful to Dani Wise for introducing us to this problem, and for discussions regarding the relevance of the MSQT. We are grateful to Vincent Guirardel for comments which greatly helped improve the exposition of the article. Thank you Sam Fisher for many enthusiastic conversations during the completion of this work. 
\end{ack}

\section{Background}

\subsection{Combinatorial complexes and words}

A map between CW-complexes is \emph{combinatorial} if $n$-cells are sent homeomorphically to $n$-cells for every $n$. A $n$-dimensional CW-complex $X$ is \emph{combinatorial} if $n=0$ or if $X^{(n-1)}$ is combinatorial and each attaching map $\partial B^n\to X^{(n-1)}$ is combinatorial for some suitable combinatorial CW-structure on boundaries of $n$-cells $\partial B^n$. In this paper, every complex will be a combinatorial, $2$-dimensional CW-complex, and every map between complexes will be combinatorial. In a complex, each $2$-cell has the cell structure of a polygon.

A \emph{path in $X$} is a (combinatorial) map $p:I\to X$ from $I$ a subdivided interval or a single $0$-cell to $X$. Paths $p_1:I_1\to X$ and $p_2:I_2\to X$ are considered equivalent if there exists an isomorphism $I_1\to I_2$ so that $I_1\xrightarrow{p_1}X$ and $I_1\to I_2\xrightarrow{p_2}X$ are the same map. The \emph{length} of a path $p$ is the number of 1-cells in $I$, which we denote by $|p|$. A path of length zero is \emph{trivial}. A path $p$ is \emph{closed} when the endpoints of $I$ are mapped to the same point of $X$. An \emph{oriented path} is a path $p:I\to X$ with a fixed orientation of $I$, and oriented paths are considered equivalent only if the isomorphism $I_1\to I_2$ is orientation-preserving. The \emph{inverse} of an oriented path $p:I\to X$ is the path with the opposite orientation of $I$, which we denote $p^{-1}$. For oriented paths $p_1$ and $p_2$, the concatenation $p_1p_2$ is well-defined provided the terminal vertex of $p_1$ and initial vertex of $p_2$ are equal. The concatenation of a path with its inverse yields a closed path that is homotopic in $X^{(1)}$ to a trivial path rel endpoints. A path decomposing as a concatenation $pp^{-1}$ with $p$ non-trivial is called \emph{fully backtracking}. The concatenation point in the domain is its \emph{midpoint}.

A \emph{cycle} is a map $c:S\to X$ from $S$ a subdivided circle or a single $0$-cell to $X$. Cycles $c_1:S_1\to X$ and $c_2:S_2\to X$ are considered equivalent if there exists an isomorphism $S_1\to S_2$ so that $S_1\xrightarrow{c_1}X$ and $S_1\to S_2\xrightarrow{c_2}X$ are the same map. The \emph{length} of a cycle is the number of $1$-cells in $S$, which we denote $|c|$. A cycle of length zero is \emph{trivial}. An \emph{oriented cycle} is a cycle $c:S\to X$ with a fixed orientation of $S$, and oriented cycles are considered equivalent only if the isomorphism $S_1\to S_2$ is orientation-preserving. A \emph{subpath} $p$ of a cycle $c:S\to X$ is a path factoring as $p:I\to S\xrightarrow{c} X$ where $I\to S$ is injective on the interior of $I$. 

A closed path $p:I\to X$ naturally factors through $S=I/\{\text{endpoints}\}$ and thus \emph{represents} a cycle. Conversely, given a cycle $c:S\to X$, for any surjective subpath $I\to S$, the concatenation $I\to S\to X$ is a closed path representing $c$. There are $|c|$ such closed paths. A path or cycle is \emph{immersed} if its defining map is locally injective. Otherwise, it contains a \emph{backtrack}, i.e.~a fully backtracking subpath. A \emph{maximal backtrack} is a fully backtracking subpath that is not properly contained in a fully backtracking subpath with the same midpoint. Note that a maximal backtrack can be properly contained in another maximal backtrack with a different midpoint.

Let $\{a_i, i\in I\}$ be a set of formal letters and $\{a_i^{-1}, i\in I\}$ their formal inverses. We will consider words in the alphabet $\{a_i^\pm\} = \{a_i, i\in I\} \sqcup \{a_i^{-1}, i\in I\}$. The \emph{inverse} of a word $w=a_1^{\epsilon_1}\cdots a_n^{\epsilon_n}$ with $a_1,\dots,a_n\in \{a_i,i\in I\}$ and $\epsilon_1,\dots,\epsilon_n\in \{+1,-1\}$ is  $w^{-1}\defeq a_n^{-\epsilon_n}\cdots a_1^{-\epsilon_1}$. Likewise, the \emph{conjugate} of a word $w_1$ by a word $w_2$ is defined as the concatenation $w_2w_1w_2^{-1}$. A word is \emph{reduced} if it does not contain a subword of the form $a_ia_i^{-1}$ or $a_i^{-1}a_i$. A word is \emph{cyclically reduced} if all of its cyclic permutations are reduced.

Whenever we orient and label the oriented $1$-cells of a complex $X$, we insist that if an oriented edge $e$ is labeled by $a_i$, then the reversed orientation of $e$ is labeled by $a_i^{-1}$. When $X$ is oriented and labeled in this way, it induces a labeling of every oriented path in $X$ by a word in the $\{a_i^\pm\}$, and a labeling of every oriented cycle in $X$ by a word in the $\{a_i^\pm\}$ up to cyclic permutation. Note that such paths (resp. cycles) are immersed when the corresponding words are reduced (resp. cyclically reduced). 

Let $G$ be a group given by a presentation $\langle a_i, i\in I\mid r_j, j\in J\rangle$. The \emph{presentation complex} of $G$ (corresponding to the given presentation) is the complex $P_G$ built as follows: $P_G$ has a single $0$-cell, $|I|$ $1$-cells with orientations labeled bijectively by the $\{a_i^\pm\}$, and $|J|$ $2$-cells whose attaching maps correspond to cycles labeled with the words $r_j$. By van~Kampen's theorem, $\pi_1P_G = G$. Each word in the $\{a_i^\pm\}$ defines a closed oriented path in $P_G$, and therefore an element of $G$. Moreover, each element of $G$ arises in this way. A covering $X\to P_G$ induces a labeling of oriented $1$-cells of $X$ by $\{a_i^\pm\}$. 

We will be interested in constructing covers of complexes. In this setting, the following criterion using links characterizes when a combinatorial map $X\to Y$ is a covering. 

\begin{defi}
    Let $X$ be a complex, $e$ a $1$-cell in $X$ and $x$ an endpoint of $e$. An \emph{end of $e$ at $x$} is a half-edge of $e$ incident to $x$. Note that every $1$-cell has two ends, even when its endpoints are the same. If $e$ is oriented, then the \emph{incoming end} of $e$ is the end towards which $e$ is oriented. The \emph{outgoing end} is the end away from which $e$ is oriented.
    
    
    Let $R$ be a $2$-cell in $X$ with its induced polygon cell structure. Then a \emph{corner} and \emph{side} of $R$ take their usual meaning, as with polygons. When $R$ is incident to a $0$-cell $x$, a \emph{corner of $R$ at $x$} is a corner of $R$ that is attached to $x$. Similarly, when $R$ is incident to a $1$-cell $e$, a \emph{side of $R$ at $e$} is a side of $R$ that is attached to $e$. Note that the $2$-cell $R$ can have multiple corners at a single $0$-cell or multiple sides at a single $1$-cell of $X$. For each side of $R$ at $e$, there are two naturally associated, distinct corners of $R$, one for each end of $e$. Likewise, for each corner $c$ of $R$ at $x$, there are two naturally associated, distinct sides of $R$ at (possibly equal) $1$-cells incident to $x$. The corner $c$ provides a natural choice of an end of these $1$-cells.

    Let $x$ be any $0$-cell of $X$. The \emph{link $\lk(x)$ of $x$} is a graph with a vertex for each end of a $1$-cell at $x$ and an edge for each corner of a $2$-cell at $x$. Intuitively, $\lk(x)$ can be viewed as the cell structure on an ``$\epsilon$-sphere'' centered at $x$. In fact, one can take a small neighborhood $N_\epsilon(x)$ of a $0$-cell $x$ in $X$ and a topological embedding $\lk(x)\to X$ so that $\partial N_\epsilon(x)$ is exactly the image of $\lk(x)\to X$. The closed neighborhood $N_\epsilon(x)\cup \partial N_\epsilon(x)$ is homeomorphic to the cone on $\lk(x)$. A combinatorial map $f\colon X\to Y$ between complexes maps sides to sides and corners to corners. Since $f$ is combinatorial, it maps ends of $1$-cells at $x$ to ends of $1$-cells at $f(x)$, and corners of $2$-cells at $x$ to corners of $2$-cells at $f(x)$. This defines a \emph{local map} $f_x\colon \lk(x)\to \lk(f(x))$ which is combinatorial. Note that $f$ also maps sides at a $1$-cell to sides at its image.
\end{defi}

\begin{rem}
    To avoid confusion, we will use "vertices" and "edges" to denote the cells of links (and, later on, the cells of underlying graphs in Definition~\ref{defi:graphofspaces}), but "$0$-cells" and "$1$-cells" to denote the cells of $X$.
\end{rem}

\begin{lm}[Combinatorial covering condition]
\label{lem:checkCover}
    Let $f\colon X\to Y$ be a combinatorial map between non-empty complexes with $Y$ connected. Then $f$ is a covering if and only if for every $0$-cell $x\in X$, the local map $f_x\colon\lk(x)\to \lk(f(x))$ is an isomorphism.
\end{lm}
\begin{proof}
    Assume first that $f$ is a covering. Then $f$ preserves the local incidence geometry at each $0$-cell of $X$. However, this local incidence geometry determines entirely the links of $0$-cells, thus all local maps $f_x$ are necessarily isomorphisms.

    Conversely, assume that every $f_x,\, x\in X^{(0)}$ is an isomorphism. Let $e$ be a $1$-cell of $X$ and choose an end of $e$ at some $0$-cell $x$, corresponding to a vertex $v$ of $\lk(x)$. Distinct sides at $e$ are associated to distinct corners at the chosen end, which in turn define distinct edges in $\lk(x)$, incident to $v$. By injectivity of $f_x$, these edges are mapped to distinct incident edges in $\lk(f(x))$, corresponding to corners associated with sides at the $1$-cell $f(e)$. Therefore, $f$ maps distinct sides to distinct sides.
    
    Let $e$ be a $1$-cell in $X$ and let $s$ be a side at $f(e)$, a $1$-cell of $Y$. Let $x$ be an endpoint of $e$. By surjectivity of $f_x$, an edge of $\lk(f(x))$ corresponding to $s$ is the image of some edge in $\lk(x)$. Since $x$ is an endpoint of $e$, this edge in $\lk(x)$ corresponds to a side at $e$, and $f$ maps this side to $s$. Thus $f$ induces a bijection between the sides at $e$ and at $f(e)$.
    
    Let $y$ be a point of $Y$ and $U$ a neighborhood of $y$. To prove that $f$ is a covering we must identify $f^{-1}(y)\times U\simeq f^{-1}(U)$ so that $f^{-1}(y)\times U\simeq f^{-1}(U) \overset{f}{\to} U$ is the second coordinate projection, provided $U$ has been chosen small enough. If $y\in f(X)$ is a $0$-cell of $Y$, $f^{-1}(U)$ is a union of neighborhoods of $0$-cells of $X$, that can be made disjoint and homeomorphic to cones of links if $U$ is small enough. Then by the assumption that local maps are isomorphisms, $f$ restricts to a homeomorphism from each component of $f^{-1}(U)$ to $U$, and the result holds at $y$. If $y\in f(X)$ is inside an open $2$-cell of $Y$, take $U$ small enough to be entirely contained in this $2$-cell. Then $f^{-1}(y)$ contains at most one point per $2$-cell of $X$, and $f^{-1}(U)$ decomposes as a union of open subsets lying in disjoint $2$-cells of $X$. As $f$ is homeomorphic on each open $2$-cell, the result holds at $y$ once again.

    Finally, assume $y\in f(X)$ lies inside an open $1$-cell $e$ of $Y$ that has $k$ distinct sides. Then, $y$ admits a small open neighborhood $U$ homeomorphic to $k$ Euclidean half-discs glued along their diameters. Distinct preimages of $y$ belong to distinct $1$-cells of $X$. Thus, provided $U$ has been chosen small enough, any component $\mathcal{C}$ of $f^{-1}(U)$ only contains one point of $f^{-1}(y)$ and is also homeomorphic to a gluing of half-discs. Since $f$ induces a bijection between sides at a $1$-cell and sides at its image, $f$ restricts to a homeomorphism $\mathcal{C}\to U$, and the result holds at $y$.

    The previous discussion shows that $f(X)\subseteq Y$ is open. Since $f(X)$ is a subcomplex of $Y$, it is closed. Since $X$ is non-empty and $Y$ is connected, $f$ is surjective, hence is a covering.
\end{proof}

\subsection{Graphs of groups}

The theory of graphs of groups, or Bass-Serre theory, is presented in \cite{serre:Trees,scottWall:TopologicalMethods}. We follow the perspective of Scott and Wall \cite{scottWall:TopologicalMethods}.

\begin{defi}
\label{defi:graphofspaces}
Let $\GG$ be a graph, i.e.~a $1$-dimensional complex, with oriented $1$-cells. Let $V$ denote the set of vertices (i.e.~$0$-cells) of $\GG$ and $E$ its set of edges (i.e.~$1$-cells). A \emph{graph of spaces} $X_\GG$ with \emph{underlying graph} $\GG$ is described by the following data: for each vertex $v$ a path-connected topological space $X_v$ called \emph{a vertex space}, and for each edge $e$ with endpoints $v_0$, $v_1$, a path-connected topological space $X_e$ called \emph{an edge space} together with $\pi_1$-injective \emph{attaching maps} $X_e\times \{0\}\to X_{v_0}$ and $X_e\times \{1\}\to X_{v_1}$. The \emph{realization} of $X_\GG$ is the quotient of $\left(\bigsqcup_{v\in V} X_v\right)\sqcup \left(\bigsqcup_{e\in E} X_e\times[0,1]\right)$ obtained by attaching $X_e\times [0,1]$ to $X_{v_0}$ and $X_{v_1}$ via the attaching maps, for each edge $e$ with endpoints $v_0$, $v_1$.

In our setup, each $X_v$, $v\in V$ and each $X_e\times [0,1]$, $e\in E$ will be complexes. Furthermore, $X_e\times \{0\}$ and $X_e\times \{1\}$ will be $1$-dimensional subcomplexes of $X_e\times [0,1]$ and every attaching map will be combinatorial. Our assumptions ensure that the realization of $X_\GG$ is always a complex. We will often confuse $X_\GG$ with its realization. The \emph{graph of groups} corresponding to a graph of spaces is its image under the $\pi_1$-functor. A \emph{splitting} of a group $G$ \emph{as a graph of groups} is an isomorphism $G\simeq\pi_1X_\GG$ for some graph of spaces $X_\GG$. The splitting is a \emph{multiple HNN-extension} if $\GG$ has a single vertex. The \emph{vertex groups} of the splitting are the groups $\pi_1X_v$ and the \emph{edge groups} are the groups $\pi_1X_e$. These groups embed as subgroups of $\pi_1X_\GG$.
\end{defi}

For any graph of spaces $X_\GG$ there is an associated continuous projection map $\pi:X_\GG\to \GG$, which sends each vertex space $X_v$ to $v$ and sends each product $X_e\times [0,1]$ to $e$ via the projection $X_e\times [0,1]\to [0,1] \simeq e$ . Any covering $\psi:\widehat{X_\GG}\to X_\GG$ of a graph of spaces has its own graph of spaces structure: The vertex spaces of this structure are the components of the preimages $\psi^{-1}(X_v)$ for all $v\in V$, and the edge spaces are given by the components of $\psi^{-1}(X_e\times \{\frac12 \})$ for all $e\in E$. Each vertex space (resp. edge space) of $\widehat{X_\GG}$ is a cover of the corresponding $X_v$ (resp. $X_e\simeq X_e\times \{\frac{1}{2}\}$). Letting $\widehat \GG$ be the underlying graph of $\widehat {X_\GG}$, there exists a map $\widehat \GG \to \GG$ such that the following diagram commutes:

\begin{center}
	\begin{tikzcd}
            \widehat{X_\GG} \arrow[r] \arrow[d] & \widehat \GG \arrow[d] \\
		X_\GG \arrow[r] & \GG
	\end{tikzcd}
    \end{center} 

The \emph{Bass-Serre tree} of $X_\GG$ (or of the splitting $\pi_1 X_\GG$) is the underlying graph $\mathcal{T}$ to the graph of spaces structure of the universal cover $\widetilde{X_\GG}$, together with the unique action of $\pi_1 X_\GG$ making the projection $\pi \colon \widetilde {X_\GG}\to \mathcal{T}$ equivariant. The graph $\mathcal{T}$ is a tree and none of its edges is inverted under this action. Stabilizers of vertices of $\mathcal{T}$ are all conjugates of vertex groups in $\pi_1 X_\GG$ and stabilizers of edges are all conjugates of edge groups.

\begin{rem}
\label{covergraph}
    As a partial converse, given a graph of spaces $X_\GG$, and a cover of graphs $\widehat\GG \to \GG$, there exists a cover $\widehat{X_\GG}\to X_\GG$ of the same degree with underlying graph $\widehat\GG$ making the previous diagram commute: The vertex spaces and edge spaces appearing in $\widehat{X_\GG}$ are the same that appear in $X_\GG$, repeated as many times as the degree of the cover $\widehat \GG \to \GG$.
\end{rem}

We are primarily interested in \emph{graphs of free groups with cyclic edge groups}, by which we mean graphs of groups with free vertex groups and infinite cyclic edge groups. In the case of a multiple HNN-extension, such a group admits presentations with abstract generators $\{a_i^\pm,i\in I,t_j^\pm,j\in J\}$ (where $I$ indexes a free basis of the vertex group and $J$ indexes edges of the underlying graph), where every relation is of the form $t_jwt_j^{-1} = w'$ for some $j\in J$ and $w$, $w'$ arbitrary words in the $\{a_i^\pm\}$. Here the vertex group is $\gen{a_i,i\in I}$ and the edge groups are all the $\gen{w}$, $\gen{w'}$. We will call any such presentation a \emph{standard presentation}, see Section~\ref{detecting}.

\subsection{Coverings and residual finiteness}

\begin{defi}
        Consider a map $f\colon A\to X$ and a covering $\widehat X\to X$. An \emph{elevation of $f$ to $\widehat X$} is a lift $\hat f\colon \widehat A \to \widehat X$, where $\widehat A$ is a minimal cover of $A$ where such a lift exists. In other words, the following diagram commutes, and there does not exist an intermediate cover $\widehat A\to \widehat A'\to A$ such that $\hat f$ factors through a lift $\widehat A'\to \widehat X$ of $f$.

    \begin{center}
	\begin{tikzcd}
            \widehat A \arrow{r}{\hat f} \arrow[d] & \widehat X \arrow{d} \\
		A \arrow{r}{f} & X
	\end{tikzcd}
    \end{center}

    An elevation corresponds to a maximal subgroup of $\pi_1A$ whose image under $f_*\colon \pi_1A\to \pi_1X$ lies in $\pi_1\widehat X<\pi_1 X$.
\end{defi}

\begin{defi}
    Let $G$ be a group, and let $X$ be a finite $2$-complex with $\pi_1X\simeq G$. The group $G$ is \emph{residually finite} if any of the following equivalent conditions is satisfied:

    \begin{enumerate}
        \item For every non-trivial $g\in G$ there exists a finite quotient $G\to \overline G$ so that the image $\bar g\in \overline G$ of $g$ is non-trivial.
        \item For every non-trivial $g_1,\dots, g_n\in G$ there exists a finite quotient $G\to \overline G$ so that the images $\bar g_1,\dots, \bar g_n\in \overline G$ are non-trivial.
        \item For every compact subcomplex $A\subset \widetilde X$ in the universal cover of $X$ there exists a finite degree cover $\widehat X\to X$ so that the restriction of $\widetilde X\to \widehat X$ to $A$ is an embedding.
        \item For each closed path $p$ in $X$ representing a non-trivial element in $G$, there exists a finite-degree cover $\widehat{X}\to X$ where some lift of $p$ is not closed.
        
    \end{enumerate}
    The equivalence is standard and proved as follows: $(1)\Rightarrow (2)$ because a finite intersection of finite-index normal subgroups is finite-index and normal. $(2)\Rightarrow (3)$ considering the finite set of $g_i$ such that $g_iA\cap A\neq \emptyset$. $(3)\Rightarrow (4)$ clearly, and $(4)\Rightarrow (1)$ because every finite-index subgroup contains a normal finite-index subgroup. For more details see \cite[Lemma~1.3]{scott:subgroups}.
\end{defi}

It follows from the definitions that subgroups of residually finite groups are residually finite, and virtually residually finite groups are residually finite.

\begin{lm}
\label{freeprodrf}
Let $A,B$ be two groups. The free product $A*B$ is residually finite if and only if both $A$ and $B$ are residually finite.
\end{lm}

\begin{proof}
Both $A$ and $B$ embed as subgroups of $A*B$. Hence if the latter is residually finite, so are the former. Conversely, assume $A$ and $B$ are residually finite and let $g\in A*B$ be non-trivial. Let $g = a_1b_1\cdots a_nb_n$ be a normal form, that is, $a_i\in A$, $b_i\in B$, and none of them are trivial except possibly $a_1$ or $b_n$ (see \cite[Chapter~IV, Theorem~1.2]{lyndonSchupp:combinatorialGroupTheory}). Consider finite quotients $\overline A$ and $\overline B$ such that none of the $\bar{a_i}$ (resp. $\bar{b_i}$) are trivial unless $a_i$ (resp. $b_i$) was trivial. Then $\bar g = \bar{a_1}\bar{b_1}\cdots \bar{a_n}\bar{b_n}$ is a normal form for $\bar g$ in the quotient $A *B\to \overline A*\overline B$. As a free product of finite groups, $\overline A*\overline B$ is virtually free \cite[Proposition~4]{serre:Trees} hence residually finite. Thus there exists a further quotient of $A*B$ to a finite group where the image of $g$ is non-trivial.
\end{proof}

\subsection{Small-cancellation theory}

There are many accounts of classical small-cancellation theory. Our discussion most closely follows \cite{mccammondWise:fansLadders}. There are some differences in definitions, but the definitions of ``reduced diagram'' and ``$C'(\lambda)$'' are equivalent to those given in \cite{mccammondWise:fansLadders}. For another account, see \cite[Chapter 5]{lyndonSchupp:combinatorialGroupTheory}. Throughout this section, we will assume that all attaching maps of $2$-cells are immersions.

\begin{defi}
Let $R$ be a closed $2$-cell in a complex $X$. The \emph{boundary cycle} $\partial_cR\colon S\to X$ of $R$ is the cycle in the complex defined by the attaching map of $R$. A \emph{boundary path} $\partial_p R\colon I\to X$ of $R$ is any closed path representing the boundary cycle.

Let $R_1$, $R_2$ be two (possibly equal) $2$-cells of $X$, with boundary cycles $\partial_cR_1\colon S_1\to X$, $\partial_cR_2\colon S_2\to X$, and let $p\colon I\to X$ be a non-trivial path immersed in $X$. The path $p$ is a \emph{piece in $X$} for $R_1$, $R_2$ if there exist maps $\iota_1\colon I\to S_1$, $\iota_2\colon I\to S_2$ such that $p$ factors both as $\partial_cR_1\circ \iota_1$ and $\partial_cR_2\circ \iota_2$, but no combinatorial isomorphism $S_1\to S_2$ makes the following diagram commute.

\begin{center}
    \begin{tikzcd}
        I \arrow{d}{\iota_1} \arrow{r}{\iota_2} & S_2 \arrow{d}{\partial_cR_2}\\
         S_1 \arrow[ru] \arrow{r}{\partial_cR_1} &X
    \end{tikzcd}
\end{center}

The complex $X$ \emph{satisfies $C'(\frac 1n)$} if for every $2$-cell $R$ and every piece $p:I\to X$ factoring through $\partial_c R$, $|p|<\frac1n |\partial_c R|$. If $\widehat X\to X$ is a cover and $p:I\to X$ is a piece, then any lift $\hat p:I\to \widehat X$ of $p$ is a piece in $\widehat X$. Hence covers of $C'(\frac 1n)$ complexes are $C'(\frac 1n)$. 
\end{defi}

\begin{defi}
A \emph{(disc) diagram} is a complex $D$ with an embedding into a complex $\overline{D}$ such that:
\begin{itemize}
    \item $\overline{D}$ is homeomorphic to $S^2$.
    \item $D$ is the complement in $\overline{D}$ of the interior of a $2$-cell, denoted $R_\infty$.
\end{itemize}

A path $p\colon I\to D$ is \emph{exterior} if the image of $p$ lies in the boundary of $D\subset\overline{D}$. A path is \emph{interior} if it maps open edges and interior vertices of $I$ to the interior of $D\subset \overline{D}$. Note some paths are neither interior nor exterior, and only trivial paths can be both. A $1$-cell $e$ of $D$ is \emph{interior} (resp. \emph{exterior)} if the corresponding path of length $1$ is interior (resp. exterior). A $2$-cell $R$ of $D$ is \emph{exterior} if it contains an exterior $1$-cell, and \emph{interior} otherwise. A $2$-cell $R$ is a \emph{shell} if some boundary path $\partial_p R$ is the concatenation of a (possibly trivial) interior path and a non-trivial exterior path in $D$. We call these paths the \emph{inner path} and the \emph{outer path} of $R$, respectively. An \emph{$i$-shell} is a shell whose inner path is the union of at most $i$ pieces. A \emph{spur} is a $0$-cell of degree $1$ in the boundary of $D\subset \overline{D}$.

A \emph{diagram over $X$} is a combinatorial map $\phi:D\to X$ for some diagram $D$. A \emph{boundary path} of a diagram over $X$ is the composition of $\phi$ with a boundary path of $R_\infty$ in $\overline{D}$. A diagram $\phi:D\to X$ over $X$ is \emph{reduced} if for every piece $p:I\to D$ in $D$, the composition $\phi \circ p:I\to X$ is a piece in $X$. In particular, if $D\to X$ is a reduced diagram and $X$ is a $C'(\frac 16)$ complex, then $D$ is $C'(\frac 16)$. 
\end{defi}

The following lemma was first proven by van Kampen \cite{vanKampen:lemmasTheoryGroups}. Thirty years later, it was rediscovered by Lyndon. An account of Lyndon's proof can be found in \cite[Chapter~5]{lyndonSchupp:combinatorialGroupTheory}. See also \cite[Lemma~2.17]{mccammondWise:fansLadders}.

\begin{lm}[van Kampen's Lemma]
\label{lem:vanKampen}
    Let $p$ be a nullhomotopic closed path in a complex $X$. There exists a reduced diagram $D\to X$ having $p$ as a boundary path.
\end{lm}

We will need a weak version of what is commonly known as ``Greendlinger's Lemma''. Stronger statements and proofs can be found in \cite[Theorem~9.5]{wise:richesToRAAGS} and \cite[Theorem~9.4]{mccammondWise:fansLadders}.

\begin{lm}[Greendlinger's Lemma]
\label{lem:greendlinger}
    Let $D$ be a reduced $C'(\frac16)$ diagram. Suppose $D$ is not a single $0$-cell. Then $D$ contains at least two 3-shells and/or spurs. 

    If $D$ has exactly two shells and/or spurs, then any shell of $D$ is a $1$-shell. 
\end{lm}

A connected subcomplex $A\subset X$ is \emph{convex} if every geodesic of $X$ joining points of $A$ is contained in $A$. It is easy to see that an intersection of convex subcomplexes is convex, thus connected.

\begin{cor}
\label{cor:relatorsEmbedded}
    Let $X$ be a simply-connected $C'(\frac16)$ complex. For any $2$-cell $R$, the inclusion map $\psi: R\to X$ is an embedding of a convex subcomplex, isometric on $1$-skeleta.
\end{cor}

\begin{proof}
    Suppose the claim is false, and let $p'$ be a shortest geodesic joining points of $\psi(R)$ which does not factor as $\psi\circ p$ for any $p:I\to X$. By Lemma~\ref{lem:vanKampen}, for any $p:I\to R$ such that $\psi\circ p$ has the same endpoints as $p'$, there exists a reduced diagram $\phi:D\to X$ having $(\psi\circ p)p'$ as a boundary path. Choose such $p$ and $D$ so that $D$ has a minimal number of cells. By assumption, $D$ has at least one $2$-cell. 

    The diagram $D$ cannot contain a spur. Indeed, $\phi$ cannot map a spur to an interior vertex of $p'$ since $p'$ is geodesic. If $\phi$ maps a spur to an interior vertex of $\phi\circ p$, then $p$ contains a backtrack, since $\psi$ is an immersion. Then, removing the backtrack from $p$ and the spur from $D$ reduces the number of cells in $D$, a contradiction. Finally, $\phi$ cannot map a spur to an endpoint of $p'$, otherwise $p'$ would not have been a shortest counterexample to the claim.

    Let $R'$ be a shell of $D$ with outer path $q:I\to D$. Assume, for contradiction, that there exists a subpath $r$ of $p$ such that $\psi\circ r$ is a shared subpath of $\phi\circ q$ and $\psi\circ p$ of length $\frac 16|\partial_c R'|$. Then, by the $C'(\frac 16)$ condition, $\phi(\partial R')=\psi(\partial R)$. In this case, replacing the subpath $r$ in $p$ by its complement in $\partial R$ and deleting $R'$ from $D$ decreases the number of cells in $D$, a contradiction. Therefore, $\phi\circ q$ and $\psi\circ p$ don't share a subpath of length $\frac 16|\partial_c R'|$. Likewise, $\phi\circ q$ and $p'$ cannot share a subpath of length greater than $\frac 12|\partial_c R'|$. Indeed, such a subpath of $p'$ could be replaced by its complement in $\partial_c R'$, reducing the length of $p'$.
    
    In summary, $\phi\circ q$ does not share a subpath longer than $\frac 16|\partial_c R'|$ with $\psi\circ p$ nor a subpath longer than $\frac 12|\partial_c R'|$ with $p'$. Yet, $\phi\circ q$ is a subpath of the cycle determined by $(\psi\circ p)p'$, a boundary path of $D$. Thus, $|\phi\circ q|\leq \frac 23 |\partial_c R'|$. If $R'$ were a $1$-shell of $D$, then $|\phi\circ q|>\frac 56 |\partial_c R'|$ would hold. Hence $D$ has no $1$-shells.

    By Lemma~\ref{lem:greendlinger}, $D$ has at least three $3$-shells. For at least one of these $3$-shells, with outerpath $q$, the path $\phi\circ q$ has to be contained in either $p'$ or $\psi \circ p$. Yet, $|\phi\circ q|>\frac 12 |\partial_c R'|$, contradicting the previous observations.
\end{proof}

\begin{cor}
\label{cor:shortLoopsAndPaths}
Let $M>0$ be the length of the shortest relator in a simply-connected $C'(\frac16)$ complex $X$.
\begin{enumerate}
    \item Any non-trivial, closed, immersed path $p$ in $X$ has length at least $M$.
    \item If $p$ is an immersed path with endpoints on a $2$-cell $R$ of $X$ and $|p|< M/2$, then $p$ is an embedded path in $R$.
\end{enumerate}
\end{cor}

\begin{proof}
    $(1)$ By Lemma~\ref{lem:vanKampen}, there exists a reduced diagram $D\to X$ having $p$ as a boundary path. Up to passing to a subpath of $p$, assume that $D$ has no spurs. If $D$ is a single 2-cell, then $|p|\geq M$. Otherwise, by Lemma~\ref{lem:greendlinger}, there exists two 3-shells $R$ and $R'$ of $D$. By the $C'(\frac16)$ condition, their outer paths have length at least $\frac12|\partial_c R|$, $\frac 12|\partial_c R'|$ respectively, in particular at least $\frac12M$ each. These outerpaths are subpaths of $p$ with disjoint interiors, hence $|p|\geq M$. 

    $(2)$ Suppose the statement is false and $p$ is a shortest counterexample. Let $q$ be a geodesic joining the endpoints of $p$. By Corollay~\ref{cor:relatorsEmbedded}, $q$ is a path in $R$. Since $p$ is shortest, the closed path $pq^{-1}$ is immersed, and contradicts $(1)$.
\end{proof}

The following theorem is a consequence of far-reaching work of Agol and Wise. We isolate a particular case of their work for our use later. 

\begin{thm}
\label{thm:smallCancellationResiduallyFinite}
    If $X$ is a $C'(\frac16)$ complex and $\pi_1X$ is finitely generated, then $\pi_1X$ is residually finite. 
\end{thm}

\begin{proof}
    By \cite[Theorem~1.2]{wise:cubulatingSmallCancellation}, $\pi_1X$ acts properly and cocompactly on a $\mathrm{CAT}(0)$ cube complex. It is well-known that $\pi_1 X$ is hyperbolic (for example, it has linear Dehn function by Lemma~\ref{lem:greendlinger}). Thus by \cite[Theorem~1.1]{agol:virtualHaken}, $\pi_1X$ is virtually special. In particular, $\pi_1X$ is linear. Thus $\pi_1X$ is residually finite by Malcev's Theorem \cite{malcev:isomorphicMatrixRepresentations}.
\end{proof}

\section{Unbalanced elements and Baumslag-Solitar subgroups}

\begin{defi}
\label{def:veryUnbalanced}
    Let $g\in G$ be an infinite-order element in a group. Say $g$ is \emph{unbalanced} for $h\in G$, $m,n\in \ZZ\setminus\{0\}$ if $|m|\neq |n|$ and $hg^mh^{-1}=g^n$ holds.

    Say $g$ is \emph{very unbalanced} for $h\in G$, $m,n\in \ZZ\setminus\{0\}$ if $1$, $|m|$ and $|n|$ are distinct and the following hold:
    \begin{itemize}
        \item $hg^mh^{-1} = g^n$
        \item If $n$ divides $m$, $hg^{m/n}h^{-1} \neq g$
        \item If $m$ divides $n$, $hgh^{-1}\neq g^{n/m}$.
    \end{itemize} 
\end{defi}

The definitions above are clearly invariant under conjugacy.

\begin{ex}
    Let $p,q$ be non-zero integers. The \emph{Baumslag-Solitar group} $BS(p,q)$ is defined by the presentation $\langle a,t \mid ta^pt^{-1}=a^q\rangle$. If $|p|\neq |q|$, then the element $a\in BS(p,q)$ is unbalanced. If $1$, $|p|$, and $|q|$ are distinct, then $a\in BS(p,q)$ is very unbalanced (this is a consequence of Britton's Lemma, see Lemma~\ref{lem:britton}).
    
    Note that the equation $ta^2t^{-1} = a^4$ holds both in $BS(2,4)$ and $BS(1,2)$, yet $a$ is very unbalanced only in $BS(2,4)$.
\end{ex}

\begin{rem}
\label{rem:rfwithvu}
    Wise (\cite{wise:separabilityGraphsFreeGroups}) identified the presence of unbalanced elements as an obstruction for any group to be LERF. Here we prove that very unbalanced elements are an obstruction for a group splitting as a graph of free groups to be residually finite. However, this is not true for an arbitrary group: Consider the subgroup of $GL_2(\QQ)$ generated by $a = \begin{pmatrix} 3/2 &0\\ 0 &1\end{pmatrix}$ and $b = \begin{pmatrix} 1 &1\\ 0 &1\end{pmatrix}$. As a finitely generated, linear group, it is residually finite by \cite{malcev:isomorphicMatrixRepresentations}, yet $b$ is very unbalanced for $a$, $2$, $3$. This group is in fact metabelian, splitting as $\ZZ\left[\frac{1}{6}\right]$-by-$\ZZ$.

\end{rem}

\begin{lm}
\label{unbalancedelliptic}
    Let $G$ split as a graph of groups, with Bass-Serre tree $\mathcal{T}$, and let $g\in G$ be unbalanced for $h,m,n$. Then $g$ fixes a vertex of $\mathcal{T}$. Moreover, if the vertex groups of the splitting are free, $h$ does not fix any vertex of $\mathcal{T}$ and $g^n$ fixes an edge of $\mathcal{T}$.
\end{lm}
\begin{proof}
    Let $d=gcd(m,n)$. Since $g^m$ and $g^n$ are conjugates in $G$, they act on $\mathcal{T}$ with the same translation length $\ell(g^m) = \ell(g^n)$. Then $\ell(g) = \ell(g^m)/m = \ell(g^n)/n = 0$, i.e.~$g$ fixes a vertex $v$ (see \cite[Proposition~25]{serre:Trees}).

    The element $g^n = hg^mh^{-1}$ fixes the geodesic segment $S$ joining $v$ and $hv$. Assuming $h$ fixes some vertex of $\mathcal{T}$, $h$ fixes a vertex $w$ of $S$. Then the free group $Stab(w)$ contains $h$ and $g^n$, hence contains $g^m$ and $g^d\in \gen{g^m,g^n}$. The element $g^{d}$ is unbalanced in $Stab(w)$ for $h,m/d,n/d$, hence fixes a vertex of the Cayley tree for $Stab(w)$ by the same translation length argument as before, a contradiction. Thus $h$ fixes no vertex, and $g^n$ fixes all the edges joining the distinct vertices $v$ and $hv$.
\end{proof}

\begin{lm}
\label{bsnovu}
    For every non-zero integer $p$, $BS(1,p)$ contains no very unbalanced elements.
\end{lm}
\begin{proof}
    Using the relation $tat^{-1}=a^{p}$,  every element $h\in BS(1,p) = \gen{a,t}$ can be represented as a product of the form $t^{-i}a^jt^k$ with $i,k\geq 0$ and $j\in \ZZ$. 
    


    Assume $g\in BS(1,p)$ is very unbalanced for $h,m,n$. By Lemma~\ref{unbalancedelliptic}, $g$ fixes a vertex of the Bass-Serre tree of the HNN extension defining $BS(1,p)$. Up to conjugating, assume $g$ is of the form $a^l$, $l\neq 0$. Let $i,k\geq 0$, $j\in\ZZ$ so that $h=t^{-i}a^jt^k$. Then the following equivalent equations hold:
    \[\begin{aligned}
        g^n &= hg^mh^{-1}\\
        a^{nl} &= t^{-i}a^jt^ka^{ml}t^{-k}a^{-j}t^i\\
        t^{i} a^{nl} t^{-i} &= a^jt^ka^{ml}t^{-k}a^{-j}\\
        a^{p^inl} &= a^j a^{p^k ml} a^{-j} = a^{p^kml}
    \end{aligned}\]

    Therefore $p^inl = p^kml$, since $a$ has infinite order. If $i<k$, then $n = p^{k-i}m$ is a multiple of $m$ and $hgh^{-1} = t^{-i}a^{p^kl}t^i = a^{p^{k-i}l} = g^{n/m}$ holds. Symmetrically if $i>k$, then $m$ is a multiple of $n$ and $hg^{m/n}h^{-1} = g$ holds. If $i=k$, then $m=n$. This contradicts the fact that $g$ is very unbalanced.
\end{proof}

For the next result, recall (see \cite{moldavanskii:residualProperties}) that $BS(p,q)$ is LERF if and only if $|p|=|q|$; $BS(p,q)$ is residually finite if and only if $1$, $|p|$, $|q|$ are not distinct; and both being LERF and residually finite are preserved under passing to a subgroup.

\begin{pro}
\label{vuimpliesbs}
    Let $G$ split as a graph of free groups.
    \begin{enumerate}
        \item $G$ contains an unbalanced element if and only if there exists $p,q\in\ZZ\setminus\{0\}$ with $|p|\neq |q|$ such that $BS(p,q)$ embeds in $G$. In that case, $G$ is not LERF.
        \item $G$ contains a very unbalanced element if and only if there exists $p,q\in \ZZ\setminus\{0\}$ with $1$, $|p|$ and $|q|$ distinct such that $BS(p,q)$ embeds in $G$. In that case, $G$ is not residually finite.
    \end{enumerate}
\end{pro}
\begin{proof}
    The "if" directions of both statements are clear since an unbalanced (resp. very unbalanced) element of a subgroup of $G$ is unbalanced (resp. very unbalanced) in $G$ as well.

    For the "only if" direction, let $g\in G$ be unbalanced for $h,m,n$ and let $d=gcd(m,n)$. Let $\mathcal{T}$ be the Bass-Serre tree of the splitting of $G$ and $\mathcal{T}'$ the subtree given by the union $\bigcup_{k\geq 1} Fix(g^k)$. The union is connected since all the terms contain $Fix(g)$. For every vertex $v$ of $\mathcal{T}'$, there exists $k\geq 1$ such that $v\in Fix(g^k)\subseteq Fix(g^{mk})$. Thus $hv\in Fix(hg^{mk}h^{-1}) = Fix(g^{nk})\subseteq \mathcal{T}'$. Therefore, $\mathcal{T}'$ is $\gen{g,h}$-stable.

    Moreover, for every $h'\in \gen{g,h}\cap Stab(v)$, $h'$ conjugates $g^{km^Kn^K}$ to $g^{kL}$, for some integers $K$ and $L$, when $K$ is sufficiently large. Note that $h'$, $g^{km^Kn^K}$ and $g^{kL}$ are in $Stab(v)$. Since $Stab(v)$ is free, $g^{km^Kn^K}$ and $g^{kL}$ have the same translation length in the $Stab(v)$ Cayley tree, hence $|km^Kn^K|=|kL|$. Therefore, $h'^2$ and $g^{kL}$ commute in $Stab(v)$, meaning that $h'$ is contained in the cyclic subgroup generated by a maximal root $r$ of $g^{kL}$ in $Stab(v)$. Since $\gen{g^{kL}}\subseteq \gen{g,h}\cap Stab(v) \subseteq \gen{r}$, the group $\gen{g,h}\cap Stab(v)$ is infinite cyclic. This holds for any vertex $v$ of $\mathcal{T}'$. Moreover, for an arbitrary edge $e$ of $\mathcal{T}'$, $\gen{g,h}\cap Stab(e)$ contains a non-trivial power of $g$ and is contained in a vertex stabilizer, thus $\gen{g,h}\cap Stab(e)$ is infinite cyclic as well.

    The action of $\gen{g,h}$ on $\mathcal{T}'$ has infinite cyclic vertex and edges stabilizers, and therefore provides a splitting of $\gen{g,h}$ as a Generalized Baumslag-Solitar group. Since $hg^mh^{-1} = g^n$ and $g$ has infinite order, $\gen{g,h}$ is a non-abelian quotient of $BS(m,n)$.

    By the argument above applied to $g^d$, which is unbalanced for $h, m/d, n/d$, $\gen{g^d,h}$ splits as a Generalized Baumslag-Solitar group and is a non-abelian quotient of $BS(m/d,n/d)$. Moreover, $|m/d|$ and $|n/d|$ are coprime. By a theorem of Levitt (\cite[Lemma~7.4]{levitt:quotients}), $BS(m/d,n/d)$ embeds in $\gen{g^{d},h}$, proving Assertion~(1).

    Now assume further that $g$ is very unbalanced for $h$, $m$, $n$. If $m$ is not a multiple or a divisor of $n$, then $1$, $|m/d|$ and $|n/d|$ are distinct, and the previous argument suffices to prove Assertion~(2) as well. Otherwise, say $m$ divides $n$, the other case being symmetric. By Lemma~\ref{bsnovu}, the Generalized Baumslag-Solitar group $\gen{g,h}$ is not isomorphic to $BS(1,p)$ for any $p$, i.e.~not solvable. Yet $BS(1, n/m)$ embeds in $\gen{g^m,h}\leq \gen{g,h}$ by the argument above. By another theorem of Levitt (\cite[Lemma~7.6]{levitt:quotients}), there exists a prime $q$ such that $BS(q, qn/m)$ embeds in $\gen{g,h}$, proving Assertion~(2) in that case.
\end{proof}

\begin{rem}
    A variant of the previous argument justifies that when $g\in G$ is conjugate into an edge group, the commensurator of $\gen{g}$ in $G$ splits as a Generalized Baumslag-Solitar group. By \cite{wise:separabilityGraphsFreeGroups}, the absence of unbalanced elements in $G$ is equivalent to all such commensurators being LERF (i.e.~\emph{unimodular} in Levitt's terminology). When the edge groups are cyclic, it is also equivalent to $G$ itself being LERF.

    Likewise, our main result implies that the absence of very unbalanced elements in $G$ is equivalent to all such commensurators being residually finite (i.e.~unimodular or isomorphic to some $BS(1,p)$). When the edge groups are cyclic, it is also equivalent to $G$ itself being residually finite.
\end{rem}

\begin{lm}
\label{freeprodvunbalanced}
Let $A,B$ be arbitrary groups. The free product $A*B$ contains a very unbalanced element if and only if $A$ or $B$ does. Moreover, for $|p|\neq |q|$, $BS(p,q)$ embeds in $A*B$ if and only if it embeds in $A$ or in $B$
\end{lm}

\begin{proof}
If $A$ contains $g$ very unbalanced for $h,m,n$, the equality $hg^mh^{-1}g^{-n} = 1$ and the non-equalities $hg^{m/k}h^{-1}g^{-n/k}\neq 1$ for $k>1$ dividing $gcd(m,n)$ hold in $A*B$, since $A$ embeds in $A*B$.

Conversely, if $g\in A*B$ is very unbalanced for $h,m,n$, then $g$ fixes a vertex in the Bass-Serre tree of the free product, fixing a vertex $v$, by Lemma~\ref{unbalancedelliptic}. The element $hg^mh^{-1} = g^n$ is non-trivial (as $g$ is infinite-order) and fixes both $v$ and $hv$. Since edge stabilizers are trivial, $v = hv$, proving that $g$ is very unbalanced for $h,m,n$ inside the subgroup $Stab(v)$ which is conjugate to either $A$ or $B$.

Moreover, if $BS(p,q)$ embeds in $A$ or $B$, it clearly embeds in $A*B$ through the monomorphisms $A\to A*B$, $B\to A*B$. Conversely, if $BS(p,q)=\gen{a,t}$ embeds in $A*B$, $a$ fixes a vertex $v$ in the Bass-Serre tree of the free product by Lemma~\ref{unbalancedelliptic}. Since edge stabilizers are trivial and $a$ has infinite order, $Fix(a^p) = Fix(a^q) =\{v\}$. However, $\{v\} = Fix(a^q) = Fix(ta^pt^{-1}) = tFix(a^p) = \{tv\}$, hence $\gen{a,t}\leq Stab(v)$. Up to conjugating, $BS(p,q) = \gen{a,t}$ embeds in $A$ or $B$.
\end{proof}

\begin{rem}
\label{rem:canUseHNN}
When $G$ splits as a graph of groups with cyclic edge groups, there exists a free group $F$ such that $G*F$ splits as a multiple HNN-extension of a free group $H$ along cyclic edge groups. This can be seen in the graph of spaces by wedging all the vertex spaces into a single vertex space, which is homotopically equivalent to wedging the space with a graph.

When the splitting of $G$ has finite underlying graph, $F$ can be taken of finite rank, and the splitting of $G*F$ has finite underlying graph as well. In this case, the union of edge groups generate a finite-rank subgroup of $H$, which is contained in a finite-rank free factor. Thus, $H$ decomposes as a free product $K_1*K_2$ such that $K_1$ is finite-rank and contains all the edge groups. This induces a new splitting of $G*F$ as $G'*K_2$, where $G'$ splits as a finite multiple HNN-extension of $K_1$ along cyclic edge groups. Note that $F$ and $K_2$ are residually finite and do not contain very unbalanced elements nor $BS(p,q)$ subgroups.

The results of this paper state implications between residual finiteness, absence of very unbalanced elements and absence of certain Baumslag-Solitar subgroups. By Lemma~\ref{freeprodrf} and Lemma~\ref{freeprodvunbalanced}, instead of working with some arbitrary graph of free groups with cyclic edge groups $G$ we can restrict to working with $G*F$, a multiple HNN-extension of a free group along cyclic edge groups. Moreover when the initial splitting of $G$ is finite, we can further restrict to working with $G'$, a finite multiple HNN-extension of a finite-rank free group along cyclic edge groups.
\end{rem}

\bin{
\section{Very unbalanced elements}

Recall that \cite{wise:separabilityGraphsFreeGroups} identified the presence of unbalanced elements as an obstruction for a group $G$ to be LERF. An \emph{unbalanced element} is an element $g\in G$ such that there exists $h\in G$, $n,m\in \ZZ$ with $|n|\neq |m|$ and $hg^mh^{-1} = g^n$. As we consider the weaker property of residual finiteness, we need a stronger obstruction.

\label{veryunbalanced}
\begin{defi}
In a group $G$, an element $g\in G$ is \emph{very unbalanced} if there exists $h\in G$, $m,n\in \ZZ$, with $1$, $|m|$, and $|n|$ distinct, such that the following hold:
\begin{itemize}
    \item $g$ has infinite order
    \item $hg^mh^{-1} = g^n$
    \item $hg^{m/k}h^{-1} \neq g^{n/k}$ for any $k>1$ dividing both $m$ and $n$.
\end{itemize} 
\end{defi}

It clearly follows from the definition that being very unbalanced is a conjugacy class invariant. Fix an arbitrary group $G$.

\begin{lm}
\label{unbalancedroot}
Let $g\in G$ be very unbalanced for $h,m,n$ and let $r\in G$ be a \emph{root} of $g$ in the sense that $r^L=g$ for some $L\in \ZZ$. Then there exists an integer $l\geq 1$ such that $r$ is very unbalanced for $h,lm,ln$.
\end{lm}
\begin{proof}
Let $L$ be such that $g=r^L$. Since $g$ has infinite order, $r$ has infinite order. Let $l\geq 1$ dividing $L$ be minimal satisfying $hr^{lm}h^{-1} = r^{ln}$. Clearly, $1,|lm|,|ln|$ are distinct.

Assume for some $k>1$ dividing $lm$ and $ln$ that $hr^{lm/k}h^{-1} = r^{ln/k}$ holds. Set $k' = k/gcd(l,k)$ and take both sides of the previous equation to the power $gcd(l,k)L/l$ to get $hr^{Lm/k'}h^{-1} = r^{Ln/k'}$. Substituting, we have $hg^{m/k'}h^{-1} = g^{n/k'}$. Hence $k' = 1$, $k$ divides $l$, and $l/k$ contradicts the minimality of $l$.
\end{proof}

The proof of the following lemma appears in \cite{moldavanskii:residualProperties}, though the lemma is not stated directly.

\begin{lm}
\label{commutatordies}
Let $\varphi$ be a homomorphism from $G$ to a finite group. Let $g,h\in G$, $m,n\in\ZZ$ and $d = gcd(m,n)$. If $hg^mh^{-1} = g^n$, then the element $\varphi([hg^dh^{-1},g])$ is trivial.
\end{lm}

\begin{proof}
Since $\varphi(hgh^{-1})$ and $\varphi(g)$ are conjugate, they have the same order $\omega$. By assumption, $\varphi(hgh^{-1})^m$ and $\varphi(g)^n$ are equal and thus have the same order $\omega/gcd(\omega, m)=\omega/gcd(\omega, n)$. Thus $gcd(\omega, m)=gcd(\omega, n)$, so $gcd(\omega, m)$ divides $d$. The elements $m$ and $gcd(\omega, m)$ generate the same cyclic subgroup of $\ZZ/\omega \ZZ$, so we have $\varphi(hgh^{-1})^d\in \gen{\varphi(hgh^{-1})^{gcd(\omega,m)}}=\gen{\varphi(hgh^{-1})^{m}}=\gen{\varphi(g)^n}\subseteq \gen{\varphi(g)}$, and the lemma follows.

\end{proof}

The following lemma contains standard, well-known results about free group elements.

\begin{lm}
\label{lem:freeroot}
Any non-trivial element $g$ in a free group $F$ has a \emph{maximal root} $r$ so that $g$ is a power of $r$ and $r$ is not a proper power. This root is unique up to inversion.

Moreover, given two non-trivial elements $g,h$ in $F$, the following are equivalent:
\begin{enumerate}
    \item $g$ and $h$ commute
    \item $g$ and $h$ have a common power
    \item $g$ and $h$ are powers of the same element
\end{enumerate}
\end{lm}
\begin{proof}
Consider a Cayley tree for $F$. If $g\in F$ is non-trivial, the subgroup $S_g<F$ of elements acting by translations on the axis of $g$ in the tree is infinite cyclic, as $S_g$ is a subgroup of the cyclic group of orientation-preserving automorphisms of a bi-infinite line graph. By uniqueness of the axis, any element with a power equal to $g$ stabilizes the axis of $g$. This means that the only two maximal roots of $g$ are the generators of $S_g$.

$(3)\Rightarrow (1)$ and $(3)\Rightarrow (2)$ are clear. If $(1)$ or $(2)$ is satisfied and $g,h$ are non-trivial, some power of $h$ preserves the axis for $g$ in a Cayley tree. By uniqueness of the axis of $h$, both $g$ and $h$ are powers of a maximal root of $g$.
\end{proof}

\begin{pro}
\label{rfimpliesnounbalanced}
Assume that $G$ splits as a graph of free groups. If $G$ is residually finite, $G$ contains no very unbalanced element.
\end{pro}

Here we make no finiteness assumptions on the splitting or the rank of the free groups. We make no assumptions on the edge groups of the splitting either.

\begin{proof}
Assume by contradiction that $G$ is residually finite and $g\in G$ is very unbalanced for $h,m,n$. Without loss of generality, say $1<n<|m|$. Set $d=gcd(m,n)\geq~1$.

As conjugates, the elements $g^m$ and $g^n$ have the same translation length in the Bass-Serre tree of the graph of groups, yet $|m|\neq|n|$, hence $g$ is elliptic, stabilizing a vertex $v$. Assume that $h$ acts elliptically as well. Then $h$ fixes a vertex $w$ on the (possibly trivial) geodesic from $v$ to $hv$. However, the element $g^n = hg^mh^{-1}$ fixes both $v$ and $hv$, hence it must fix $w$. Thus the free group $Stab(w)$ contains $h,g^m,g^n$, and $g^d \in \gen{g^m,g^n}$. In a Cayley tree for $Stab(w)$, the elements $g^m$ and $g^n$ are powers of $g^d$ with the same translation lengths. This is a contradiction as $g$ is of infinite order and $|n|<|m|$. Thus, moving forward, we suppose that $h$ acts loxodromically on the Bass-Serre tree. Let $\alpha_h$ be its axis, with ends denoted $h^{\pm \infty}$.

\bigskip

Assume first that $Fix(g)$ does not contain the end $h^{+\infty}$, and let $x$ be the vertex of $Fix(g)$ closest to $h^{+\infty}$. More precisely:
\begin{itemize}
    \item If $Fix(g)$ doesn't intersect $\alpha_h$, let $x$ be the vertex of $Fix(g)$ closest to $\alpha_h$. See Figure~\ref{fig:FixTrees1}.
    \item Otherwise, $Fix(g)$ intersects $\alpha_h$ along an interval that does not have $h^{+\infty}$ as an endpoint, let $x$ be the endpoint of that interval closest to $h^{+\infty}$. See Figure~\ref{fig:FixTrees2}.
\end{itemize}

Let $r$ be a maximal root of $g$ in the free group $Stab(x)$, so that $x\in Fix(r)$. Since $Fix(r)\subseteq Fix(g)$, $Fix(r)$ does not contain $h^{+\infty}$, and $x$ is the point of $Fix(r)$ closest to $h^{+\infty}$. Observe that, in particular, the geodesic segment $\gamma$ joining $x$ to $hx$ only intersects $Fix(r)$ at $x$.

\begin{figure}[h]
    \centering
    \begin{subfigure}[b]{0.45\textwidth}
        \includegraphics[width=\textwidth]{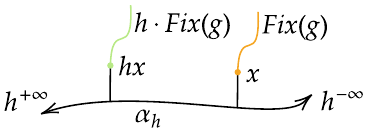}
        \caption{}
        \label{fig:FixTrees1}
    \end{subfigure}
    \begin{subfigure}[b]{0.45\textwidth}
        \includegraphics[width=\textwidth]{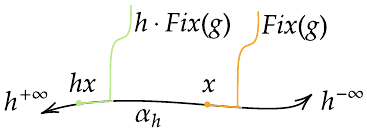}
        \caption{}
        \label{fig:FixTrees2}
    \end{subfigure}
    \caption{Different cases for the position of the point $x$, depending on whether $Fix(g)$ intersects $\alpha_h$ or not.}
\end{figure}

By Lemma~\ref{unbalancedroot}, $r$ is very unbalanced for $h$ and some $lm, ln$. By Lemma~\ref{commutatordies} and residual finiteness of $G$, $[hr^{ld}h^{-1},r]$ is trivial in $G$. Hence $Fix(r)$ and $Fix(hr^{ld}h^{-1})$ intersect at some vertex $y$. The set $Fix(hr^{ld}h^{-1})$ is convex and contains both $hx$ and $y$. Thus it contains the geodesic segment $\gamma'$ joining these two points. Since $y\in Fix(r)$, $\gamma'$ contains $\gamma$. This proves that $x$ belongs to $Fix(hr^{ld}h^{-1})$. Therefore, $r$ and $hr^{ld}h^{-1}$ both lie in the free group $Stab(x)$.

Since $r$ and $hr^{ld}h^{-1}$ are commuting elements of a free group and $r$ is not a proper power, there exists $p\in\ZZ$ such that $hr^{ld}h^{-1} = r^p$ by Lemma~\ref{lem:freeroot}. Then $r^{ln} = (hr^{ld}h^{-1})^{m/d} = r^{pm/d}$ and $ln = pm/d$, since $r$ has infinite order. Setting $k = m/d = ln/p$, the equation $hr^{lm/k}h^{-1} = r^{ln/k}$ holds, and because $r$ is very unbalanced, $|k|=1$ and thus $d=|m|$, a contradiction with the assumption that $1<n<|m|$.
\medskip

Assume now that $Fix(g)$ contains the end $h^{+\infty}$. Let $x$ be a vertex in $\alpha_h\cap Fix(g)$, and $r$ a maximal root of $g$ in $Stab(x)$. As before, $r$ is very unbalanced for $h$ and some $lm,ln$. If $Fix(r)$ does not contain the end $h^{+\infty}$ the previous argument applied to $r$ instead of $g$ yields a contradiction. Otherwise, $Fix(r)$ contains $hx$, hence $x\in Fix(h^{-1}rh)\cap Fix(r^{ld})$. As before $[r,hr^{ld}h^{-1}]$ is trivial in $G$, and so is its conjugate $[h^{-1}rh,r^{ld}]$. By Lemma~\ref{lem:freeroot} applied to $Stab(x)$ and the fact that $r$ is not a proper power in $Stab(x)$, there exists $p\in \ZZ$ such that $h^{-1}rh = r^p$. Thus $r^{pln} = h^{-1}r^{ln}h = r^{lm}$ and $pln = lm$ because $r$ has infinite order. Dividing, we see that $lm/ln=p$ is an integer. The equation $hr^{lm/ln}h^{-1} = r^{ln/ln}$ and the fact that $r$ is very unbalanced imply that $|ln| = 1$, a contradiction.
\end{proof}

\begin{rem}
\label{rem:resfinwithveryunbalanced}
    Note that without the splitting assumption, a residually finite group can indeed contain a very unbalanced element: Consider the subgroup of $GL_2(\QQ)$ generated by $a = \begin{pmatrix} 2/3 &0\\ 0 &1\end{pmatrix}$ and $b = \begin{pmatrix} 1 &1\\ 0 &1\end{pmatrix}$. As a finitely generated, linear group, it is residually finite by \cite{malcev:isomorphicMatrixRepresentations}, yet $a$ conjugates $b^2$ and $b^3$. This group is in fact metabelian, splitting as $\ZZ\left[\frac{1}{6}\right]$-by-$\ZZ$.
\end{rem}

\begin{lm}
\label{freeprodvunbalanced}
Let $A,B$ be arbitrary groups. The free product $A*B$ contains a very unbalanced element if and only if $A$ or $B$ does.
\end{lm}

\begin{proof}
If $A$ contains $g$ very unbalanced for $h,m,n$, the equality $hg^mh^{-1}g^{-n} = 1$ and the non-equalities $hg^{m/k}h^{-1}g^{-n/k}\neq 1$ for $k>1$ dividing $gcd(m,n)$ hold in $A*B$, since $A$ embeds in $A*B$.

Conversely, if $g\in A*B$ is very unbalanced for $h,m,n$, then $g$ is elliptic in the Bass-Serre tree of the free product, fixing a vertex $v$, as in the proof of Proposition~\ref{rfimpliesnounbalanced}. The element $hg^mh^{-1} = g^n$ is non-trivial (as $g$ is infinite-order) and fixes both $v$ and $hv$. Since edge stabilizers are trivial, $v = hv$, proving that $g$ is very unbalanced for $h,m,n$ inside the subgroup $Stab(v)$ which is conjugate to either $A$ or $B$.
\end{proof}

\begin{rem}
\label{rem:canUseHNN}
When $G$ splits as a graph of groups with cyclic edge groups, there exists a free group $F$ such that $G*F$ splits as a multiple HNN-extension of a free group $H$ along cyclic edge groups. This can be seen in the graph of spaces by wedging all the vertex spaces into a single vertex space, which is homotopically equivalent to wedging the space with a graph.

When the splitting of $G$ has finite underlying graph, $F$ can be taken of finite rank, and the splitting of $G*F$ has finite underlying graph as well. In this case, the union of edge groups generate a finite-rank subgroup of $H$, which is contained in a finite-rank free factor. Thus, $H$ decomposes as a free product $K_1*K_2$ such that $K_1$ is finite-rank and contains all the edge groups. This induces a new splitting of $G*F$ as $G'*K_2$, where $G'$ splits as a finite multiple HNN-extension of $K_1$ along cyclic edge groups. Note that $F$ and $K_2$ are residually finite and do not contain very unbalanced elements (by Proposition~\ref{rfimpliesnounbalanced}).

The results of this paper state implications between residual finiteness and absence of very unbalanced elements. By Lemma~\ref{freeprodrf} and Lemma~\ref{freeprodvunbalanced}, instead of working with some arbitrary graph of free groups with cyclic edge groups $G$ we can restrict to working with $G*F$, a multiple HNN-extension of a free group along cyclic edge groups. Moreover when the initial splitting of $G$ is finite, we can further restrict to working with $G'$, a finite multiple HNN-extension of a finite-rank free group along cyclic edge groups.
\end{rem}
}

\section{Detecting very unbalanced elements}
\label{detecting}

\begin{defi}
Let $G$ split as a multiple HNN-extension of a free group $F$ with cyclic edge groups. Let $(a_i)_{i\in I}$ be a free basis of $F$. A \emph{standard presentation} of $G$ is of the form $G=\langle (a_i)_{i\in I},\, (t_j)_{j\in J}\mid t_j u_j^{m_j} t_j^{-1} = v_j^{n_j}\text{ for }j\in J\rangle$, where $J$ is a set, the $m_j$, $n_j$ are integers, and the $u_j$, $v_j$ are (not necessarily distinct) words in the $\{a_i^\pm\}$ representing elements of $F$ that are not proper powers in $F$. The $t_j$ are the \emph{stable letters} of the presentation.
\end{defi}

Note that we do not assume in this section that $I$ and $J$ are finite.

\begin{defi}
In a free group $F$, a set $S$ of non-trivial group elements is \emph{independent} (or a \emph{malnormal collection}) when no distinct elements of $S$ have conjugate powers (except for the trivial power, i.e.~the identity).
\end{defi}

\begin{rem}
\label{alignhnn}
Let $j\in J$ and let $w$ be any word in the $\{a_i^\pm\}$. Replacing the generator $t_j$ with $t'_j=t_jw^{-1}$, the relation involving $t_j$ becomes $t'_j(wu_jw^{-1})^{m_j}t_j'^{-1} = v_j^{n_j}$. Consequently, $u_j$ can be replaced by an arbitrary conjugate word in the presentation without changing the resulting group. The word $u_j$ can also be replaced by its inverse up to reversing the sign of $m_j$. Symmetrically, the transformation $t'_j = wt_j$ allows one to replace $v_j$ by an arbitrary conjugate, or by its inverse up to reversing the sign of $n_j$.

Let $A\subset F$ be minimal so that for each $g\in F$ represented by a $u_j$ or $v_j$, either $g$ or $g^{-1}$ is conjugate to a (unique) element of $A$. Clearly, $|A|\leq 2|J|$ and since the group elements represented by the $u_j$ and $v_j$ are not proper powers in $F$, $A$ is independent.
For $g\in A$, let $w_g$ denote a cyclically reduced word in the $\{a_i^\pm\}$ representing $g$. By the previous observation, the presentation can be modified so that every $u_j$ and every $v_j$ is equal to some $w_g$, $g\in A$. From now on, we will assume standard presentations are of this form. In particular, each standard presentation will come with a set $A$ as above.
\end{rem}

The following lemma provides a sufficient condition for a word in the standard presentation to represent a non-trivial element of $G$ (see for example \cite[Chapter~IV, Section~2, p.181]{lyndonSchupp:combinatorialGroupTheory}).

\begin{lm}[Britton's Lemma]
\label{lem:britton}
Let $w$ be a non-empty reduced word in the $\{a_i^\pm, t_j^\pm\}$. Assume that for every subword of $w$ of the form $t_ju_j^{k}t_j^{-1}$ (resp. $t_j^{-1}v_j^{k}t_j$), $m_j$ (resp. $n_j$) does not divide $k$. Then $w$ represents a non-trivial element of $G$.
\end{lm}

\begin{defi}[The graph $\Gamma$]
\label{gammaGraph}
With the conventions of Remark~\ref{alignhnn}, given a standard presentation of $G$, define an oriented (multi-)graph $\Gamma$ with vertex set $A$ and edge set $J$, where each relation $t_j {w_g}^{m_j} t_j^{-1} = w_h^{n_j}$ provides an oriented edge $j$ from vertex $h$ to vertex $g$. Ends of edges admit a canonical labeling by non-zero integers: The terminal end of $j$ is labeled $m_j$ and the initial end is labeled $n_j$. A warning: paths and cycles in $\Gamma$ will never be required to respect the orientation of $\Gamma$, even when they are themselves oriented. 

An oriented path $\overset{\to}{p}$ in $\Gamma$ defines a word in the $\{t_j^\pm\}$. The group element \emph{represented by $\overset{\to}{p}$} is the element of $G$ represented by this word. The labels on ends of edges in $\Gamma$ induces labelings on the ends of edges of paths and cycles in $\Gamma$. Let $\overset{\to}{\gamma}$ be an oriented cycle in $\Gamma$. The \emph{loop-product} $lp(\overset{\to}{\gamma})\in \ZZ\setminus\{0\}$ is the product of all the labels on the outgoing ends of edges of $\overset{\to}{\gamma}$. A cycle $\gamma$ is \emph{balanced} if $|lp(\overset{\to}{\gamma})| = |lp(\overset{\leftarrow}{\gamma})|$, where $\overset{\to}{\gamma}$ and $\overset{\leftarrow}{\gamma}$ represent the two possible orientations of $\gamma$, and \emph{unbalanced} otherwise. A connected component $C$ of $\Gamma$ is \emph{clean} if all the cycles in $C$ are balanced. Note that it suffices to check cleanliness on cycles embedded in $C$. The terminology is different from Wise's clean graph of graphs \cite[Definition~4.5]{wise:separabilityGraphsFreeGroups} but related in the following sense: if every component of $\Gamma$ is finite and clean, the splitting satisfies the assumptions of the main theorem of \cite{wise:separabilityGraphsFreeGroups}.
\end{defi}

\begin{rem}
    In fact, each component of $\Gamma$ is the defining graph of a Generalized Baumslag-Solitar group, which embeds in $G$ as the commensurator of any vertex of the component. Every non-cyclic commensurator of an element of $G$ appears in this way. Clean components correspond exactly to the cases where the Generalized Baumslag-Solitar subgroup is \emph{unimodular} in Levitt's terminology (\cite{levitt:quotients}), i.e.~LERF.
\end{rem}

\begin{defi}
\label{potential}
    Let $C$ be a finite clean component of $\Gamma$. A \emph{potential} on $C$ is an assignment of a positive integer $p_g$ to every vertex $g$ of $C$ such that the following holds: If an edge joins $g_1$ to $g_2$ and has ends labeled $m_1$ and $m_2$ respectively, then $p_{g_1}/|m_1|$ and $p_{g_2}/|m_2|$ are equal and are integers. The existence of a potential is an easy consequence of the cleanliness and finiteness of $C$. We will henceforth endow every finite clean component of $\Gamma$ with a potential.
\end{defi}

\begin{lm}
\label{nonclean1}\begin{enumerate}
\item Let $\gamma$ be an unbalanced cycle embedded in $\Gamma$. If both loop-products of $\gamma$ are different from $\pm 1$, some vertex of $\gamma$ is a very unbalanced element of $G$.
\item Let $p$ be a closed path immersed in $\Gamma$, and $\gamma$ the corresponding cycle. If $\gamma$ is unbalanced and both ends of $p$ have labels different from $\pm 1$, the endpoint of $p$ is a very unbalanced element of $G$.
\end{enumerate}
\end{lm}

\begin{proof}
The proofs of the two statements will proceed in similar ways. We begin by choosing vertices $g$ and $h$ of $\Gamma$ and two oriented paths $\overset{\to}{p}$ and $\overset{\to}{q}$ from $g$ to $h$ in both set-ups.  

In the setting of the first statement, let $g$ and $h$ be vertices of $\gamma$ that cut $\gamma$ into two embedded oriented paths $\overset{\to}{p}$ and $\overset{\to}{q}$, both from $g$ to $h$ with the following properties: The ends of $\overset{\to}{p}$ are not labeled with $\pm 1$, and $\overset{\to}{q}$ is either trivial or all of its edges have both ends labeled with $\pm 1$. In setting of the second statement, let $g=h$ be the endpoint of $p$, let $\overset{\to}{p}$ be an arbitrary orientation of $p$, and let $\overset{\to}{q}$ be the trivial path at $g$.

\medskip

Let $k_p$ and $k_q$ be the elements of $G$ represented by $\overset{\to}{p}$ and $\overset{\to}{q}$, let $n^+$ be the loop-product of $\gamma$ with the orientation given by $\overset{\to}{p}$, and let $n^-$ be the loop-product of $\gamma$ with the opposite orientation. By assumption, the integers $1$, $|n^-|$, and $|n^+|$ are distinct. By definition of $\Gamma$, there exists $\epsilon\in \{-1,1\}$ such that the relations $g^{n^+} = k_ph^{\epsilon n^-}k_p^{-1}$ and $g = k_qh^\epsilon k_q^{-1}$ hold, yielding $k_pk_q^{-1}g^{n^-}k_qk_p^{-1} = g^{n^+}$.

Let $k\geq 1$ dividing $n^-$ and $n^+$ be maximal so that $k_pk_q^{-1}g^{n^-/k}k_qk_p^{-1} = g^{n^+/k}$ holds. If $k\neq |n^-|,|n^+|$, then $g$ is very unbalanced for $k_pk_q^{-1}$, $n^-/k$ and $n^+/k$. If $k = |n^-|$, then $k_phk_p^{-1}g^{-\epsilon n^+/n^-}=id_G$ holds. However, $\overset{\to}{p}$ is immersed and its ends are not labeled by $\pm 1$. By Lemma~\ref{lem:britton}, the left-hand side cannot be trivial, yielding a contradiction. The same argument holds for $k = |n^+|$, in which case we deduce $k_qk_p^{-1} g k_pk_q^{-1}g^{-n^-/n^+}=id_G$ and Lemma~\ref{lem:britton} provides the same contradiction.
\end{proof}

Recall that for an oriented edge $e$ in a graph, the outgoing end of $e$ is the half-edge away from which $e$ is oriented.

\begin{pro}
\label{nonclean2}
Assume $G$ contains no very unbalanced element, and let $C$ be a non-clean component of $\Gamma$. Then $C$ contains a unique embedded cycle $\gamma$, and $\gamma$ has one loop-product equal to $\pm 1$. Let $\overset{\to}{p}$ be an immersed, oriented path in $C$ meeting $\gamma$ only at its terminal vertex. Then the outgoing ends of all edges of $\overset{\to}{p}$ have labels $\pm 1$. Consequently, any path immersed in $C$ has one end labeled $\pm 1$.
\end{pro}

\begin{proof}
Fix an unbalanced cycle $\gamma$ embedded in $C$. By Lemma~\ref{nonclean1}, one of the loop-products of $\gamma$ is $\pm 1$. Since $\gamma$ is unbalanced, there exists an edge $e$ of $\gamma$ with an end not labeled by $\pm 1$. Let $x$ be the corresponding endpoint. Suppose for contradiction that there exists a cycle $\gamma'\neq \gamma$ embedded in $C$. Two cases arise: 

\begin{itemize}
    \item If $\gamma$ and $\gamma'$ intersect, up to replacing $\gamma'$ by some other embedded cycle in the union of images of $\gamma$ and $\gamma'$, assume that $\gamma'$ does not contain $e$. Let $q$ be a path starting at $x$ and traveling through $e$, then along $\gamma$ until meeting $\gamma'$ for the first time, then looping around $\gamma'$ once (in either orientation), and going back to $x$ via the same edges as before.
    \item If $\gamma$ and $\gamma'$ don't intersect, let $r$ be a path of minimal length from the image of $\gamma$ to the image of $\gamma'$. Let $y$ be the starting vertex of $r$ on $\gamma$. Let $q$ be a path starting at $x$, traveling through $e$ then along $\gamma$ until $y$, then traversing $r$, looping around $\gamma'$ once (in either orientation), and going back to $x$ via the same edges as before.
\end{itemize}

In both cases, both ends of $q$ have a label different from $\pm 1$, and $q$ is immersed. Let $s$ be the concatenation of the closed path based at $x$ looping around $\gamma$ starting with the edge $e$, and the path $q$. The path $s$ is closed, immersed, and has end labels different from $\pm 1$. Since $\gamma$ is unbalanced, at least one of $\{q, s\}$ corresponds to an unbalanced cycle. By Lemma~\ref{nonclean1}, this gives rise to a very unbalanced element, a contradiction. Hence $\gamma$ is the only embedded cycle in $C$.

Now, let $\overset{\to}{p}$ be as in the statement, and consider two closed paths $\overset{\to}{p}\gamma\overset{\to}{p}^{-1}$ and $\overset{\to}{p}\gamma^2\overset{\to}{p}^{-1}$ (with an arbitrary orientation of $\gamma$). As before, at least one of these paths corresponds to an unbalanced cycle. By Lemma~\ref{nonclean1}, the labels of ends of this path cannot both be different from $\pm 1$. Since the two ends are equal to the starting end of $\overset{\to}{p}$, this proves our claim for the first edge of $\overset{\to}{p}$. The general claim follows from the fact that every edge of $\overset{\to}{p}$ is the first edge of a subpath of $\overset{\to}{p}$ also meeting $\gamma$ at only at its terminal vertex.

Finally, since $\gamma$ is the only embedded cycle of $C$, any path immersed in $C$ either is immersed in the image of $\gamma$ or shares an end with an immersed path meeting $\gamma$ only at its other end. Thus one end of such a path is labeled $\pm 1$ (using that $\gamma$ has a loop-product equal to $\pm 1$ for the first case).
\end{proof}

\bin{\begin{cor}
    Let $G$ be a generalized Baumslag-Solitar group, i.e.~a group splitting as a finite, non-empty, connected graph of infinite cyclic groups with cyclic edge groups. Then $G$ is residually finite if and only if it is virtually $F_n\times \ZZ$ for some $n\geq 0$ or isomorphic to $BS(1,m)$ for some $m\in \ZZ\setminus\{0\}$.
\end{cor}
\com{See Levitt Subgroups and quotients of BS gps (2015) p.3}
The statement of the above corollary was remarked by Cornulier in a MathOverflow post \cite{cornulier:overflow}.
\begin{proof}
    $F_n\times \ZZ$ is a direct product of residually finite groups and thus is residually finite. $BS(1,m)$ is residually finite as well \cite{baumslagSolitar:theOGpaper, moldavanskii:residualProperties}. See also Example~\ref{ex:BS(1,q)}.  
    
    Suppose that $G$ is a residually finite generalized Baumslag-Solitar group. Either $G\simeq \ZZ$, in which case the statement is obvious, or $G$ is the fundamental group $\pi_1 X_{\GG}$ of a graph of spaces $X_{\GG}$ with vertex and edge spaces homeomorphic to $S^1$, such that the attaching maps are covering maps, and $\GG$ is not a single vertex. Using Remark~\ref{rem:canUseHNN}, define a group $\widehat G = G*F$ for some free group $F$, splitting as a multiple HNN-extension of a free group. Let $\Gamma$ be the graph associated to that multiple HNN-extension. Since $G$ is residually finite, $\widehat G$ contains no very unbalanced elements by Proposition~\ref{rfimpliesnounbalanced}. Hence every non-clean component of $\Gamma$ is of the form given by Proposition~\ref{nonclean2}.

    In the connected graph of spaces $X_{\GG}$ corresponding to the initial splitting of $G$, every vertex space of $X_{\GG}$ has an incident edge space. Hence, given a standard presentation of $\widehat G$ and the associated graph $\Gamma$, the vertex set $A$ of $\Gamma$ corresponds naturally to the vertex set of $\GG$. More precisely, if an element $g\in A$ corresponds to a vertex $v$ of $\mathcal{G}$, then $g\in G$ and $g$ is a generator of the vertex group at $v$. Furthermore, two vertices of $\Gamma$ are adjacent if and only if the corresponding vertex spaces are joined by an edge space in $X_{\GG}$. Thus $\Gamma$ is naturally isomorphic to $\GG$, in particular, $\Gamma$ is connected.

    First suppose $\Gamma$ is non-clean. Then by Proposition~\ref{nonclean2}, $\Gamma$ contains a unique embedded, oriented cycle $\overset{\to}{\gamma}$ with outgoing loop-product $\pm 1$ and incoming loop-product $m$. Every vertex $g$ outside the image of $\gamma$ is joined to that image by a unique immersed, oriented path $\overset{\to}{p}$ meeting $\gamma$ only at its terminal vertex, and the outgoing ends of all edges of $\overset{\to}{p}$ have labels $\pm 1$. Taking $g$ at maximal distance from the image of $\gamma$, $g$ is of degree $1$ by uniqueness of $\overset{\to}{p}$. The unique edge at $g$ corresponds to a cylinder in $X_{\GG}$, glued along the cycle corresponding to $g$ by a homeomorphism. This cycle is a free face of that cylinder. Collapse that cylinder with a homotopy equivalence onto its other boundary component. Proceeding inductively, $X_{\GG}$ collapses to its subspace induced by the image of $\gamma$, seen as a cycle in the underlying graph $\GG$. This subspace is made of cylinders glued with a homeomorphism on one boundary component and a covering map on the other. By van Kampen's theorem, the fundamental group of the subspace is isomorphic to $BS(\pm 1,m)\simeq BS(1,\pm m)$, while being isomorphic to $G$, by the homotopy equivalence. This proves the result in the first case.

    Now suppose $\Gamma$ is clean. In the terminology of Wise, the graph corresponding to $X_{\GG}$ is balanced (see \cite[Definition~4.16]{wise:separabilityGraphsFreeGroups}). By \cite[Definition~4.5, Theorem~4.18]{wise:separabilityGraphsFreeGroups}, there exists a finite-degree clean cover $Y$ of $X_\GG$. The cover $Y$ splits as a graph of spaces whose vertex spaces and edge spaces are finite-degree covers of those of $X_\GG$, thus still homeomorphic to $S^1$. Up to taking a further degree $2$ cover (to prevent immersed non-orientable surfaces), $Y$ is homeomorphic to the direct product of a graph with $S^1$, so $\pi_1 Y=F_n\times \ZZ$ is finite-index in $G$.
\end{proof}}

\begin{rem}
\label{rem:algorithm}
The rest of the paper will prove a converse of Proposition~\ref{nonclean2} for a multiple HNN-extension with finitely many edges. More explicitly, we will prove that whenever the graph $\Gamma$ satisfies the conclusion of Proposition~\ref{nonclean2} and $J$ is finite, then $G$ is residually finite, thus $G$ cannot contain very unbalanced elements by Proposition~\ref{vuimpliesbs}. This provides a simple algorithm for deciding residual finiteness of $G$ given a standard presentation with $J$ finite: first modify the presentation as in Remark~\ref{alignhnn} using the solution to the conjugacy problem in free groups, and build the graph $\Gamma$. Identify the clean and non-clean components of $\Gamma$. If every non-clean component has the structure from the conclusion of Proposition~\ref{nonclean2}, then $G$ is residually finite, otherwise it is not.
\end{rem}

\bin{\begin{rem}
\todo{apparently this explicitly contradicts a result of Levitt, Subgroups and Quotients of GBS (2015), lemma 7.4}
Very unbalanced elements for integers $n,m$ provide natural maps from $BS(n,m)$ into $G$, but these maps need not be injective, as seen in the following example. Consider the group 
\[G=\gen{a,b,s,t\mid sa^2s^{-1} = b^3, tb^5t^{-1} = a^7}\]

The graph $\Gamma$ corresponding to this standard presentation is a length $2$ cycle whose two loop-products are $10$ and $21$ (see Figure~\ref{fig:gamma_example}), thus $G$ contains very unbalanced elements. Note that $G = H*\ZZ$ where $H$ is a generalized Baumslag-Solitar group with underlying graph $\Gamma$ and vertex groups $\gen{a}$ and $\gen{b}$. 

\begin{center}
\begin{figure}[h]
    \begin{tikzcd}[column sep = 8em]
        \gen{a} \ar[r,bend left]{}{\text{\large{$s$}}}[pos = 0.1]{2}[pos = 0.9]{3} &\ar[l,bend left]{}{\text{\large{$t$}}}[pos = 0.1]{5}[pos = 0.9]{7} \gen{b}
    \end{tikzcd}
    \caption{The graph $\Gamma$ for $G$, coinciding with the underlying graph for the splitting of $H$.}
    \label{fig:gamma_example}
\end{figure}
\end{center}

Assume $G$ has a subgroup isomorphic to any $BS(n,m) = \gen{x,y\mid yx^ny^{-1} = x^m}$ with $1,|n|,|m|$ distinct, in the form of an injective map $f\colon BS(n,m) \hookrightarrow G$. The equation $f(y)f(x)^nf(y)^{-1}=f(x)^m$ holds. However, for $k>1$ dividing $n$ and $m$, $yx^{n/k}y^{-1}x^{-m/k}$ is non-trivial in $BS(n,m)$ by Lemma~\ref{lem:britton}. Thus by injectivity of $f$, $f(y)f(x)^{n/k}f(y)^{-1}\neq f(x)^{m/k}$. Hence $f(x)$ is very unbalanced for $f(y)$, $n$, $m$.

As in the proof of Proposition~\ref{rfimpliesnounbalanced}, $f(x)$ is elliptic in $G$ and even in a conjugate of $H$ by Remark~\ref{freeprodvunbalanced}. This implies that up to composing $f$ with a conjugation in $G$, we can assume that $f(x)\in \gen{a}$ or $\gen{b}$. Likewise, $f(y)$ is loxodromic, and up to changing the presentation of $G$\todo{this is where the problem is}, we can assume that $f(y)$ is a power of $ts$ or $st$ depending on whether $f(x)\in \gen{a}$ or $\gen{b}$, respectively. In the first case, the relation becomes, for some non-zero integers $\alpha,\beta$:
\[(ts)^\beta a^{\alpha n}(ts)^{-\beta}a^{-\alpha m} = id_G\]

Assuming $\beta>0$, $\alpha n$ is even by Lemma~\ref{lem:britton}. The relation can be rewritten:

\[(ts)^{\beta-1}t b^{3\alpha n/2}t^{-1}(ts)^{1-\beta}a^{-\alpha m} = id_G\]

By Lemma~\ref{lem:britton} again, $3\alpha n/2$ is a multiple of $5$ and the relation becomes: 
\[(ts)^{\beta-1}a^{21\alpha n/10}(ts)^{1-\beta}a^{-\alpha m} = id_G\]

Iterating the process finally yields $a^{(21/10)^\beta \alpha n - \alpha m} = id_G$, which implies $21^\beta n = 10^\beta m$ since $G$ is torsion-free. A similar equality is obtained if $\beta<0$, or if $f(x)\in \gen{b}$. This proves that $n$ is a multiple of $10^\beta$ and $m$ a multiple of $21^\beta$.

However, the following computation holds in $G$: \[\begin{aligned}\relax[(ts)^\beta (a^\alpha)^{2^\beta 5^{\beta-1}}(ts)^{-\beta}, (a^\alpha)^7] &= [t(b^{\alpha})^{3^\beta 7^{\beta-1}}t^{-1},(a^\alpha)^7] \\
&= t[(b^{\alpha})^{3^\beta 7^{\beta-1}}, (b^\alpha)^5]t^{-1}\\
&= id_G
\end{aligned}\]

Injectivity of $f$ implies that $[y x^{2^\beta 5^{\beta-1}} y^{-1}, x^7] = id$ holds in $BS(n,m)$. However, because $n$ and $m$ are multiples of $10^\beta$ and $21^\beta$ respectively, Lemma~\ref{lem:britton} yields a contradiction. Thus this particular group $G$ is not residually finite, but all of its Baumslag-Solitar subgroups are.
\end{rem}}

\section{Residual finiteness}

We start with an instructive proof of residual finiteness for $BS(1,q)$.

\begin{ex}
\label{ex:BS(1,q)}
Let $X$ be the presentation complex for $BS(1,q)$, with $1$-cells labeled by $a$ and $t$. The $1$-cells are both loops. The $2$-cell $R$ of $X$ has boundary path $tat^{-1}a^{-q}$. We build a parametrized family of covers $\widehat X(n,m)$ for $n,m>0$ with $n$ coprime to $q$. Given such $n,m$, let $\alpha \colon C\to X$ be the oriented cycle in $X$ wrapping $n$ times around the $1$-cell labeled by $a$. Let $(C_i)_{i\in \ZZ/ m\ZZ}$ be $m$ copies of $C$ and $\alpha_i$ corresponding copies of $\alpha$. Build a \emph{tube} $T$, homeomorphic to an annulus, from $n$ copies of $R$ glued along their $1$-cells labeled $t$, so that the two boundary components are labeled by $a^n$ and $a^{nq}$. Let $(T_i)_{i\in\ZZ/m\ZZ}$ be $m$ copies of $T$. For $i\in \ZZ/m\ZZ$, attach $T_i$ to $C_i$ and $C_{i+1}$, using a degree $1$ cover for $C_i$ and a degree $np$ cover for $C_{i+1}$. The resulting space is $\widehat X(n,m)$. See Figure~\ref{fig:BScover}. The map $\bigsqcup \alpha_i$ with domain $\bigsqcup C_i\subset \widehat X(n,m)$ extends to a label-preserving covering map $\widehat X(n,m)\to X$.

\begin{figure}[h]
    \centering
    \includegraphics[width=0.5\textwidth]{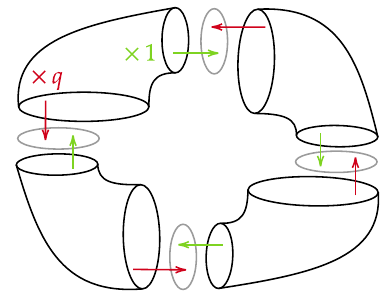}
    \caption{An example of a $\widehat X(n,m)$ cover for $BS(1,q)$.}
    \label{fig:BScover}
\end{figure}

Let $p$ be a path in $X(n,m)$ labeled by a word $w$ in $a$, $t$ representing a non-trivial element of $BS(1,q)$. Let $\hat p$ be a lift of $p$ to $\widehat X(n,m)$. If the sum of powers of $t$ in $w$ is non-zero, then for $m$ sufficiently large, $\hat p$ is a non-closed path in $\widehat X(n,m)$. If the sum of powers of $t$ in $w$ is zero, then $w$ represents a conjugate of a power of $a$ (see the proof of Lemma~\ref{bsnovu}). 
Thus for $n$ sufficiently large, $\hat p$ is not closed in $\widehat X(n,m)$. Since $w$ represented an arbitrary non-trivial element of $BS(1,q)$, this proves residual finiteness of $BS(1,q)$.
\end{ex}
\bigskip

In this section, let $G$ be a group with no very unbalanced elements splitting as a finite multiple HNN-extension with finite-rank free vertex group $F$ and cyclic edge groups. Endow $G$ with a standard presentation as in Section~\ref{detecting}. Note that in the present case, $I$ and $J$ are finite, thus $G$ is finitely presented. The goal of this section is to prove the following proposition:

\begin{pro}
\label{mainprop}
Let $w$ be a (freely) reduced word in the $\{a_i^\pm, t_j^\pm\}$ representing a non-trivial element in $G$. Then there exists a finite-degree cover of the presentation complex $P_G$ in which some path labeled by $w$ is not closed.
\end{pro}

Residual finiteness of $G$ will follow directly.

\subsection{Normalizing subwords}

\begin{defi}
Let $T$ be the total number of occurrences of stable letters $t_j^\pm$ in $w$. Up to changing $w$, assume this number is minimal among all reduced words representing the same group element.

Let $F$, $A$, and $\Gamma$ be as in Section~\ref{detecting}, and let $\mathcal{B}$ be the underlying graph of the presentation complex $P_G$, which is a bouquet of $|J|$ circles.
\end{defi}

\begin{defi}
Let $p$ be a non-trivial path in $\Gamma$, traversing the vertices $g_0,\dots, g_{l}\in A$ and edges $j_1,\dots, j_l \in J$ in order. The path $p$ represents a word $t_{j_1}^{\epsilon_1}\cdots t_{j_l}^{\epsilon_l}$ where $\epsilon_1,\dots, \epsilon_l\in \{\pm 1\}$. A \emph{$\Gamma$-word corresponding to $p$} is a reduced word of the form $w_{g_0}^{n_0}t_{j_1}^{\epsilon_1}\cdots t_{j_l}^{\epsilon_l}w_{g_l}^{n_l}$, where $n_0,\dots, n_l\in \ZZ$ (geometrically, this corresponds to a path in the Bass-Serre tree whose edges have commensurable stabilizers). A $\Gamma$-word is \emph{clean} when the image of $p$ lies in a clean component of $\Gamma$, and \emph{non-clean} otherwise.

An \emph{extremal subword} of $w$ is a subword $s$ that decomposes as the concatenation of two non-clean $\Gamma$-words $s_1$, $s_2$ whose corresponding non-trivial paths $p_1$, $p_2$ in $\Gamma$ are exactly inverse to each other, and so that $s$ is not properly contained in any subword with a decomposition of the same kind. Note that since $w$ is reduced, $s$ is itself a $\Gamma$-word.

A \emph{vertex excursion} is a maximal, non-trivial subword $u$ of $w$ that does not contain any stable letters. The \emph{apex} of an extremal subword $s$ of $w$ is the vertex excursion appearing in $s$ between the last stable letter of $s_1$ and the first stable letter of $s_2$. The apex is a power of some $w_g$, $g\in A$.
\end{defi}

\begin{figure}[h]
    \centering
    \includegraphics[width=0.8\textwidth]{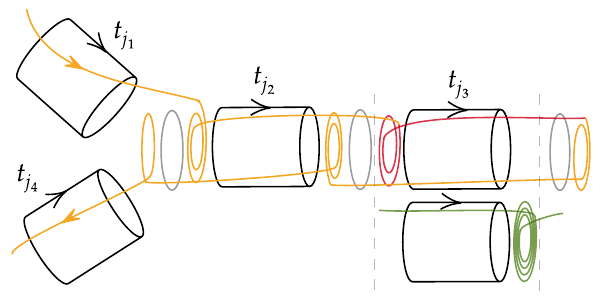}
    \caption{A visual representation of one step of the normalization procedure, replacing the red subpath with the green.}
    \label{fig:normalize}
\end{figure}

\begin{ex}
The graph $\Gamma$ corresponding to $BS(1,4)=\gen{a,t \mid tat^{-1}=a^4}$ is an oriented loop with outgoing and incoming ends labeled by 1 and 4. Consider the word $w=t^3at^{-1}a^{-3}t^{-1}a^5t^2$. The subwords $s_1=at^{-1}a^{-3}t^{-1}a^2$ and $s_2=a^3t^2$ are $\Gamma$-words corresponding to paths $p_1$ and $p_2$. The path $p_1$ wraps around $\Gamma$ twice in the direction of orientation, and $p_2$ wraps around twice in the opposite direction. The subwords $s_1$ and $s_2$ are chosen 
to be maximal so that $p_1$ and $p_2$ satisfy $p_1=p_2^{-1}$. Thus $s=s_1s_2$ is extremal. The apex of $s$ is $a^5$. Note that $s'_1 = at^{-1}a^{-3}t^{-1}a^5$ and $s'_2 = t^2$ is also a valid decomposition.
\end{ex}

\begin{lm}
\label{gammabacktrack}
Let $s$ be a $\Gamma$-word corresponding to a path $p$ in $\Gamma$ and appearing as a subword of $w$. Then for every backtrack $qq^{-1}$ of $p$, the terminal end of $q$ is not labeled by $\pm 1$. If $p$ lies in a non-clean component of $\Gamma$ and the label of the terminal end of $p$ is not $\pm 1$, then $p$ is immersed and the outgoing end of each edge of $p$ is labeled by $\pm 1$.
\end{lm}
\begin{proof}
    Let $t_j^\epsilon w_g^n t_j^{-\epsilon}$ be the smallest subword of $s$ that projects to the concatenation of the last edge of $q$ and the first edge of $q^{-1}$, where $\epsilon\in \{\pm 1\}$, $g\in A$, $n\in \ZZ$. Let $k$ and $l$ be the labels of the incoming and outgoing ends of the last edge of $q$ respectively. Then a relation $t_j^\epsilon w_g^k t_j^{-\epsilon} = w_h^{l}$ holds in $G$ for some $h\in A$. Yet $k$ does not divide $n$, otherwise substituting $t_j^\epsilon w_g^n t_j^{-\epsilon} = w_h^{nl/k}$ in $s$ would produce a word representing the same group element as $w$ with two fewer stable letters. In particular, $k\neq \pm 1$.

    Assume that $p$ is non-clean and the label of the terminal end of $p$ is not $\pm 1$. Let $p'$ be the largest immersed subpath of $p$ containing the last edge. If $p'$ is properly contained in $p$, $p$ backtracks at the first vertex of $p'$, and the label of the initial end of $p'$ is not $\pm 1$. In that case, path $p'$ is then immersed and has no end labeled $\pm 1$. Since $G$ has no very unbalanced elements by assumption, this contradicts Proposition~\ref{nonclean2}. Therefore $p=p'$ is immersed. By Proposition~\ref{nonclean2} again, the initial end of every subpath $p''$ of $p$ containing the last edge is labeled $\pm 1$, i.e. the outgoing end of each edge of $p$.
\end{proof}

\begin{lm}
\label{normalize}
Let $s = s_1s_2$ be a concatenation of two non-clean $\Gamma$-words whose corresponding paths $p_1$, $p_2$ in $\Gamma$ are inverses (in particular, $s$ is a non-clean $\Gamma$-word itself). Then there exists a reduced word $\hat s$ representing the same group element as $s$ so that $\hat s$ is a conjugate of a reduced word over $\{a_i^\pm\}$ by the word over $\{t_j^\pm\}$ labeling $p_1$. We call $\hat s$ a \emph{normalized form} of $s$.
\end{lm}

\begin{proof}
Let $k$ be the label of the terminal end of $p_1$. By Lemma~\ref{gammabacktrack} applied to $s$, we have $k\notin \{\pm 1\}$, and by the same lemma applied to $s_1$, the outgoing end of every edge of $p_1$ is labeled $\pm 1$.

Consider any subword of $s_1$ of the form $w_g^nt_j^\epsilon$ for some $j$, $g$, $n$, $\epsilon$. The edge of $p_1$ corresponding to the stable letter in this subword has outgoing end labeled $\pm 1$, thus the relation $w_g t_j^{\epsilon} = t_j^{\epsilon} w_h^{l}$ holds in $G$ for some $h\in A, l\in \ZZ$. Consequently, $s_1$ can be modified by replacing the subword $w_g^nt_j^\epsilon$ with $t_j^\epsilon w_h^{nl}$ and canceling, so that the occurrence of $t_j$ is closer from the start of $s_1$. A symmetric modification using $t_j^{\epsilon}w_g  =  w_h^{l}t_j^{\epsilon}$ can be performed in $s_2$. Performing a series of such modifications to group all the stable letters of $s_1$ as a prefix and the stable letters of $s_2$ as a suffix terminates. Deleting then pairs of adjacent inverse letters in the new apex until it is reduced yields a normalized form for $s$. Note that the number of stable letters in $s$ did not change.
\end{proof}

\begin{rem}
Let $s$ be an extremal subword of $w$. Replacing $s$ with a normalized form $\hat s$ in $w$ obtained from Lemma~\ref{normalize} produces a word $\hat w$ corresponding to the same path as $w$ in the underlying graph $\mathcal{B}$ of $P_G$, and representing the same group element as $w$. We call this operation \emph{normalizing $s$ in $w$}. Additionally, $\hat w$ is still reduced, by minimality of the number of stable letters in $w$, and normalizing $s$ did not introduce new stable letters. Thus, normalizing an extremal subword of $w$ does not affect the running assumptions on $w$. If $p$ is a path in $P_G$ with label $w$, then the normalization procedure can be visualized as pushing certain subpaths of $p$ through cylinders labeled by stable letters appearing in an extremal subword. See Figure~\ref{fig:normalize}.
\end{rem}

\begin{cor}
\label{cor:extremalnormalized}
There exists a reduced word $\hat w$ representing the same conjugacy class as $w$ such that $w$ and $\hat w$ correspond to the same path in the underlying graph $\mathcal{B}$, and all extremal subwords of $\hat w$ are normalized.
\end{cor}

\begin{proof}
First we claim that if $s$ is an extremal subword of $w$, the apex of $s$ is disjoint from any other extremal subword of $w$. Indeed, assume the apex of $s$ overlaps $s'$ and let $p_1$, $p_2$, $p'_1$, $p'_2$ be the paths in $\Gamma$ corresponding to $s$ and $s'$ respectively. Three cases arise:

\begin{itemize}
    \item If the apex of $s$ overlaps the apex of $s'$, then the apexes are equal by definition of an apex. Thus, the union of $s$ and $s'$ is itself an extremal subword. By maximality of extremal subwords, $s=s'$.
    \item If the apex of $s$ appears before the apex of $s'$ in $s'$, then some non-trivial prefix $q$ of $p_2$ is a suffix of $p'_1$. Since $s$ and $s'$ are non-clean $\Gamma$-words, by Lemma~\ref{gammabacktrack} none of the ends of $q$ can be labeled by $\pm 1$, but $q$ corresponds to a $\Gamma$-word, contradicting the same lemma.
    \item If the apex of $s$ appears after the apex of $s'$ in $s'$, then some non-trivial suffix of $p_1$ is a prefix of $p'_2$. A symmetric argument yields the same contradiction.
\end{itemize}  

As the three cases are exhaustive, the claim is proved.

We prove the corollary inductively on the number of vertex excursions in $w$, then on the sum of the lengths of vertex excursions which are not apexes.

Assume $w$ contains an extremal subword $s$ that is not normalized. If $s$ contains an entire vertex excursion that is different from the apex, normalizing $s$ decreases the total number of vertex excursions and the induction proceeds. If $s$ does not contain an entire vertex excursion besides the apex, $s$ contains a non-trivial prefix or suffix $u$ using only the $\{a_i^\pm\}$, otherwise it would already be normalized. By the previous claim, the vertex excursion containing $u$ cannot be an apex, thus normalizing $s$ does not increase the number of vertex excursions and decreases the sum of lengths of vertex excursions which are not apexes, and the induction proceeds. When no further inductive step can be executed, every extremal subword of the obtained word $\hat w$ is normalized, proving the corollary.
\end{proof}

Replacing $w$ with $\hat w$ does not affect our running assumptions on $w$: We will henceforth assume that all extremal subwords of $w$ are normalized.

\subsection{Building the cover}

\begin{center}
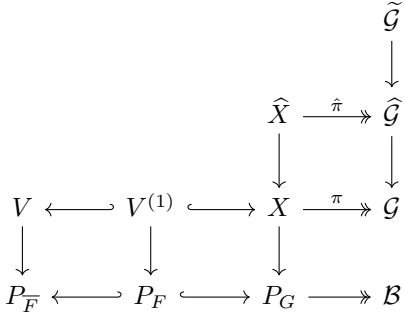
\begin{figure}[h]
	\begin{tikzcd}
        & & & \widetilde \GG \arrow{d} \\
        & & \widehat X \arrow[twoheadrightarrow]{r}{\hat \pi} \arrow{d} & \widehat \GG \arrow{d}\\
        V \arrow[hookleftarrow]{r} \arrow{d} & V^{(1)} \arrow[hookrightarrow]{r} \arrow{d} & X \arrow{d} \arrow[twoheadrightarrow]{r}{\pi} & \GG\\
        P_{\overline F} \arrow[hookleftarrow]{r} & P_F \arrow[hookrightarrow]{r} & P_G \arrow[twoheadrightarrow]{r} & \mathcal{B}
	\end{tikzcd}
 \caption{The various spaces used in the construction of the cover $\widehat X\to P_G$. Squares commute and vertical arrows are coverings.}
 \label{fig:bigDiagram}
\end{figure}
\end{center}

Recall that $P_G$, $P_F$ denote the presentation $2$-complexes for $G$ and $F$, respectively, with their standard presentations. $P_G$ has a single $0$-cell, $|I| + |J|$ $1$-cells, and $|J|$ $2$-cells, and $P_F$ is a bouquet of $|I|$ circles. Recall that $F$ is a finite-rank free group, $A\subset F$ is a finite independent subset, and $(w_g)_{g\in A}$ is a family of cyclically reduced words in the $\{a_i^\pm\}$, such that $w_g$ represents an element of $G$ conjugate to $g$.

\begin{lm}
\label{lem:boundedPiece}
There exists a constant $C>0$ so that the following holds: If for some $g,h\in A$, some integers $k$, $l$, and some cyclic permutations $w'_{g}$, $w'_{h}$ of $w_{g}$, $w_{h}$ the words $(w_{g}')^{k}$ and $(w_{h}')^{l}$ share a prefix of length at least $C$, then $g = h$ and $w'_{g}=w'_{h}$.
\end{lm}
 
\begin{proof}
Assume cyclic permutations $w'_{g}$ and $w'_{h}$ of $w_{g}$ and $w_{h}$ have powers that share the same prefix of length $lcm(|w'_{g}|,|w'_{h}|)$. The prefix represents both a power of $w'_{g}$ and a power of $w'_{h}$. Hence $w_{g}$ and $w_{h}$ have conjugate powers and $g = h$ by independence. Consequently, $|w'_{g}|=|w'_{h}|$, hence $w'_{g}=w'_{h}$. Thus it suffices to take $C=\max_{g,h\in A}\{lcm(|w'_{g}|,|w'_{h}|)\}$.
\end{proof}

\begin{lm}
\label{lem:MSQT}
Let $C>0$ be the constant from Lemma~\ref{lem:boundedPiece}. For all $g\in A$, let $n_g$ be any integer satisfying $|w_g^{n_g}|>6C$, and consider the following group presentation $\overline F := \gen{a_i\mid w_g^{n_g}}$. The corresponding presentation complex $P_{\overline F}$ is $C'(\frac16)$. In particular, $\overline F$ is a residually finite quotient of $F$.
\end{lm}

This lemma relates to \emph{omnipotence} of free groups, but in a stronger form (see \cite[Theorem~3.5]{wise:separabilityGraphsFreeGroups}). Note that it is stated here in the greatest generality for the exponents $n_g$, however we will choose these exponents more carefully in Definition~\ref{def:scquotient}.

\begin{proof}
    Let $p:I\to P_{\overline F}$ be an immersed path factoring through the boundary cycles of two $2$-cells $R$ and $R'$ of $P_{\overline F}$. Let $w_g^{n_g}$ and $w_{h}^{n_{h}}$ be labels of $\partial_c R$ and $\partial_c R'$, respectively. The cycle $\partial_c R$ in $P^{(1)}_{\overline F}$ elevates to an immersion of a line $\RR\to \widetilde{P^{(1)}_{\overline F}}$ labeled by $w_g^\infty$. Similarly, $\partial_c R'$ elevates to a line labeled by $w_{h}^\infty$, and $p$ lifts to a segment lying in the intersection of these two lines. Thus, the label of $p$ is a subword of sufficiently large powers of $w_g$ and $w_h$.

    If $|p|\geq  C$, then by Lemma~\ref{lem:boundedPiece}, $g=h$ and $w_g=w_{h}$. This implies the attaching maps $\partial_c R$ and $\partial_c R'$ are equal, and $p$ is not a piece.
    Since every $2$-cell has boundary path longer than $6C$ and every piece is shorter than $C$, the complex $P_{\overline F}$ is $C'(\frac16)$. Residual finiteness of $\overline{F}$ follows from Theorem~\ref{thm:smallCancellationResiduallyFinite}.
\end{proof}

\begin{defi}
\label{def:scquotient}
Recall $T$ is the number of occurrences of stable letters $t_j^\pm$ in $w$, and $C$ is the constant obtained from Lemma~\ref{lem:boundedPiece}. Let $N$ be an upper bound for all labels of $\Gamma$ in absolute value and $L$ be an upper bound for the length of all $w_g$, $g\in A$. Let $W=\max(|w|,2)$. Let $K=2W+3C+1$.

Since $\Gamma$ is finite, every vertex $g$ in a clean component of $\Gamma$ has a potential $p_g$ (Definition~\ref{potential}), which is a non-zero integer. For each $g\in A$, choose an integer $n_g$ satisfying the following conditions:

\begin{itemize}
    \item For every clean component $\mathcal C$ of $\Gamma$, all the vertices $g$ of $\mathcal{C}$ satisfy $n_g = (2W+KLN^T+6C+1)|p_g|$.
    \item For every non-clean component $\mathcal{C}$ of $\Gamma$, all the vertices $g$ of $\mathcal{C}$ have the same value $n_g$. Moreover, $n_g$ is greater than $2W+KLN^T+6C$ and coprime to every end label appearing in $\mathcal{C}$.
\end{itemize}

Let $\overline F=\gen{a_i \mid w_g^{n_g}}$ as in Lemma~\ref{lem:MSQT}, and $P_{\overline F}$ the corresponding presentation complex, with a single $0$-cell, $|I|$ $1$-cells and $|J|$ $2$-cells. By Lemma~\ref{lem:MSQT}, $P_{\overline F}$ is $C'(\frac16)$, and $\overline F$ is residually finite.

Let $O$ be the maximum of all $n_g$, $g\in A$, and let $\mathcal{W}$ be the finite set of words in the $\{a_i^\pm\}$ of length at most $(W+1)(LO+1)+KLN^T$ that represent non-trivial elements of $\overline F$. By residual finiteness of $\overline F$, there exists a cover $V\to P_{\overline F}$ with finite degree $d$, such that every path in $V$ with label in $\mathcal W$ is non-closed. The restriction of the covering to the $1$-skeleton $V^{(1)}$ is a cover of the graph $P_F = (P_{\overline{F}})^{(1)}$. See Figure~\ref{fig:bigDiagram}.
\end{defi}

\begin{rem}
\label{rem:shortLoopsLiftToLoops}
    Any non-closed path $\gamma$ in the universal cover $\widetilde{P_{\overline F}}$ such that $|\gamma| \leq (W+1)(LO+1)+KLN^T$ is labeled by a word in $\mathcal{W}$. By choice of $V$, the projection of $\gamma$ under $\widetilde{P_{\overline F}}\to V$ is a non-closed path. Therefore, closed paths in $V$ of length at most $(W+1)(LO+1)+KLN^T$ must lift to closed paths in the universal cover $\widetilde{P_{\overline F}}$.
\end{rem}

\begin{defi}
    Let $g\in A$. An \emph{elevation of $w_g$} to some cover $Y$ of $P_F$ is a oriented cycle elevating to $Y$ the oriented cycle labeled $w_g$ in $P_F$. In other words, it is an oriented cycle in $Y$ for which some representative closed path projects to a path in $P_F$ labeled with a power of $w_g$, but no closed subpath has the same property.
\end{defi}

\begin{lm}
\label{lem:elevationsEmbed}
For any $g\in A$, any elevation of $w_g$ to $V^{(1)}$ is an embedded cycle of length $n_g|w_g|$. Moreover, $w_g$ has exactly $d/n_g$ distinct elevations.
\end{lm}

\begin{proof}
There exists an elevation of the boundary cycle of the $2$-cell labeled $w_g^{n_g}$ at every $0$-cell of $\widetilde {P_{\overline F}}^{(1)}$. By Corollary~\ref{cor:relatorsEmbedded}, such an elevation is an embedded cycle labeled $w_g^{n_g}$. Thus every elevation of $w_g$ to $\widetilde {P_{\overline F}}^{(1)}$ is an embedded cycle labeled $w_g^{n_g}$, hence of length at most $LO$. By Remark~\ref{rem:shortLoopsLiftToLoops}, elevations of $w_g$ to $V^{(1)}$ cannot be shorter than $|w_g^{n_g}|$. Thus elevations of $w_g$ to $V^{(1)}$ are all labeled $w_g^{n_g}$ and embedded, again by Remark~\ref{rem:shortLoopsLiftToLoops}. 

The $2$-cell of $P_{\overline{F}}$ with boundary label $w_g^{n_g}$ has $d$ lifts in $V$, and every path labeled $w_g^{n_g}$ in $V^{(1)}$ arises as the boundary of one of those lifts. Yet, since $w_g$ is not a proper power, exactly $n_g$ lifts of the $2$-cell share the same boundary path, giving the desired count of elevations.
\end{proof}

\begin{lm}
\label{lem:shortinelevation}
    Let $p$ be a path immersed in $V$ of length at most $W$. Assume there exists $g\in A$ and $e_g$ some elevation of $w_g$ to $V^{(1)}$ whose image contains both endpoints of $p$. Then $p$ is embedded and $p$ or $p^{-1}$ is a subpath of $e_g$.
\end{lm}
\begin{proof}
    By Lemma~\ref{lem:elevationsEmbed}, elevations of $w_g$ to $V^{(1)}\subset V$ are boundary cycles of embedded $2$-cells labeled by $w_g^{n_g}$. Let $R$ be such a $2$-cell with boundary cycle $e_g$, so that $p$ is an immersed path with endpoints on $R$. Let $q$ be a path embedded in the boundary of $R$ joining the endpoints of $p$. Then the closed path $pq^{-1}$ has length at most $W+ LO$. By Remark~\ref{rem:shortLoopsLiftToLoops}, it lifts to a closed path $\tilde p\tilde q^{-1}$ in the universal cover $\widetilde V=\widetilde {P_{\overline F}}$. The path $\tilde p$ has endpoints on the boundary of the lift $\widetilde R$ of $R$ containing $\tilde q$. The universal cover $\widetilde V$ is a simply-connected $C'(\frac16)$ complex by Lemma~\ref{lem:MSQT}, and $|\tilde p|\leq W < n_g/2$ for all $g$, in particular $\tilde p$ is shorter than half of the boundary length of any relator. Hence, by Corollary~\ref{cor:shortLoopsAndPaths} (2), $\tilde p$ is embedded in the boundary of $\widetilde R$. Projecting, $p$ is embedded in the boundary of $R$. 
\end{proof}

The following technical result will be useful in the proof of the proposition.

\begin{lm}
\label{lem:butterflyBalloon}
Let $2\leq k \leq W$ be an integer. Let $(x_i)_{1\leq i\leq k+1}, (y_i)_{0\leq i\leq k}$ be two families of $0$-cells of $V$. Let $(g_i)_{1\leq i\leq k}$ be a family of elements of $A$ and for each $1\leq i\leq k$, $e_i$ some elevation of $w_{g_i}$ to $V^{(1)}$. Assume the following ($\dist$ denotes the combinatorial distance in $V^{(1)}$):
\begin{enumerate}
    \item For $1\leq i\leq k$, $x_i$ and $y_i$ belong to the image of $e_i$.
    \item For $1\leq i\leq k$, every path $p_i$ embedded in the image of $e_i$, joining $x_i$ to $y_i$ is longer than $W+3C$.
    \item $\displaystyle \sum_{i=0}^{k} \dist(y_i,x_{i+1})\leq W$
\end{enumerate}
and either
\begin{enumerate}
    \bin{\item[4a.] $\dist(y_0,x_{k+1})\leq W$}
    \item[(4a)] $y_0 = x_{k+1}$; or
    \item[(4b)] There exists $g\in A$ and $e_g$ an elevation of $w_g$ to $V^{(1)}$ whose image contains $y_0$ and $x_{k+1}$.
\end{enumerate}
Then for some $1\leq i \leq k-1$, $e_i=e_{i+1}$. In particular, $y_i$ and $x_{i+1}$ belong to the image of the same elevation. 
\end{lm}

\begin{figure}[h]
    \centering
    \includegraphics[width=0.8\textwidth]{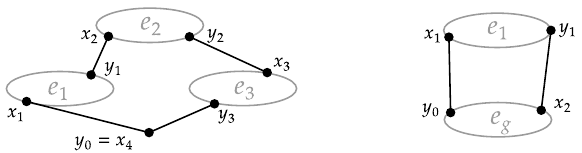}
    \caption{Two situations in which Lemma~\ref{lem:butterflyBalloon} is applicable.}
\end{figure}

\begin{proof}
Assume for the sake of contradiction that $e_i\neq e_{i+1}$ for every $1\leq i\leq k-1$. For every $1\leq i\leq k$, choose a path $p_i$ embedded in the image of $e_i$ from $x_i$ to $y_i$. Note there are two choices for each $p_i$. For every $0\leq i\leq k$, let $q_i$ be a geodesic in $V^{(1)}$ from $y_i$ to $x_{i+1}$. Finally, choose a path $r$ from $x_{k+1}$ to $y_0$ that is either trivial under Assumption~4a or embedded in the image of $e_g$ under Assumption~4b, in which case there are two choices for $r$. By Assumptions~2 and 3, $q_0p_1q_1\cdots p_{k}q_{k}r$ is a closed path of length at most $(W+1)LO+W$ in $V$. By Remark~\ref{rem:shortLoopsLiftToLoops}, its lift $\widetilde{q_0}\widetilde{p_1}\widetilde{q_1}\cdots \widetilde{p_k}\widetilde{q_k}\widetilde{r}$ to the $C'(\frac 16)$ complex $\widetilde{P_{\overline F}}$ is a closed path. 

For every $1\leq i\leq k$, $|p_i|-|q_{i-1}|-|q_i|\geq |p_i|-W> 3C>0$. Canceling all backtracks involving edges of the $q_i$, there exists a path $q'_0p'_1q'_1\cdots p'_{k}q'_{k}$ homotopic rel. endpoints to $\widetilde{q_0}\widetilde{p_1}\widetilde{q_1}\cdots \widetilde{p_{k}}\widetilde{q_{k}}$, where $q_i'$ is a subpath of $\widetilde{q_i}$ and $p_i'$ is a subpath of $\widetilde{p_i}$ with the following properties: for each $i$, $|p_i'|\geq 3C$, and if $q'_i$ is non-trivial, $p'_iq'_ip'_{i+1}$ does not backtrack. When $q'_i$ is trivial, $p'_iq'_ip'_{i+1} = p'_ip'_{i+1}$ possibly backtracks. However, this backtrack is not longer than the overlap between $e_i$ and $e_{i+1}$. Since $e_i\neq e_{i+1}$ by assumption, this overlap is no longer than $C$ by Lemma~\ref{lem:boundedPiece}. Cancelling all remaining backtracks, $q'_0p'_1q'_1\cdots p'_{k}q'_{k}$ is homotopic rel. endpoints to a non-backtracking, i.e.~immersed, path $\gamma = q''_0p''_1q''_1\cdots p''_{k}q''_{k}$ where $p''_i$ is a subpath of $p'_i$ with $|p''_i|\geq |p'_i| - 2C >C$, and $q''_i=q'_i$ unless both of them are trivial paths. Note that $\gamma$ is the unique immersed path homotopic in the $1$-skeleton of $\widetilde{P_{\overline F}}$, rel. endpoints to $\widetilde{q_0}\widetilde{p_1}\widetilde{q_1}\cdots \widetilde{p_k}\widetilde{q_k}$.

By Lemma~\ref{lem:vanKampen}, there exists $D$ a reduced diagram over $\widetilde{P_{\overline F}}$ with boundary path $\gamma r$. Choose such a $D$ of minimal area across all choices of the $p_i$, $r$, and $D$ itself. By Lemma~\ref{lem:greendlinger}, three cases arise, and we obtain contradictions in each case.

\begin{itemize}
    \item Assume $D$ has zero area, then $\gamma r$ is nullhomotopic in the $1$-skeleton of $\widetilde{P_{\overline F}}$. Since both $\gamma$ and $r$ are immersed paths, they must be inverses. But then $r$ shares subpaths of length greater than $C$ with each $p''_i$. Under Assumption~4a, this is impossible as $|r| = 0$. Under Assumption~4b, this is impossible because $e_1\neq e_2$ cannot be both equal to $e_g$, contradicting Lemma~\ref{lem:boundedPiece}.
 
    \item Assume $D$ contains a $3$-shell $R$. Under the map $D\to \widetilde{P_{\overline F}}$, the boundary of $R$ maps to an elevation $e_h$ of some $w_h$ for $h\in A$. By definition of $n_h$, $|e_h|>6C+W$, hence the outer path of the $3$-shell $R$ is longer than $3C+W$. Under Assumption~4a, $e_h$ overlaps the cycle defined by $\gamma$ along a subpath longer than $3C+W$. Since $\sum |q''_i|\leq \sum |q_i|\leq W$, $e_h$ overlaps some $p''_i$ along a subpath longer than $C$, and $e_h = e_i$ by Lemma~\ref{lem:boundedPiece}. Likewise, under Assumption~4b,  $e_h$ overlaps the cycle defined by $\gamma r$ along a subpath of length greater than $3C+W$, hence $e_h$ overlaps some $p''_i$ or $r$ along a subpath of length greater than $C$, and $e_h=e_i$ or $e_g$ by Lemma~\ref{lem:boundedPiece}. Indeed, if $|r|>C$, the argument in the previous case applies, and if $|r|\leq C$, $|r|+\sum |q''_i|\leq W+C$ is sufficient to conclude as before. In either case, changing the choice of $p_i$ and/or $r$ allows us to delete $R$ from $D$, reducing its area, a contradiction.

    \item Assume $D$ is a single $2$-cell $R$. As above, the boundary of $R$ maps to an elevation $e_h$ of some $w_h$ for $h\in A$. Under either Assumption~4a or 4b, $p''_1$ is a subpath of length greater than $C$ of $e_h$, hence $e_h = e_1$. Thus, changing the choice of $p_1$ allows us to delete $R$ from $D$, reducing its area, a contradiction.
\end{itemize}

\end{proof}

\begin{defi}
\label{defi:tubeConstruction}
We will now imitate Example~\ref{ex:BS(1,q)} by attaching tubes to $V^{(1)}$. The tubes are made of copies of $2$-cells appearing in $P_G$. For any edge $j$ in $\Gamma$ corresponding to a relator $t_jw_{g}^k t_j^{-1}w_{h}^{-l}$, two possibilities arise:

\begin{figure}[h]
    \centering
    \includegraphics[width=0.8\textwidth]{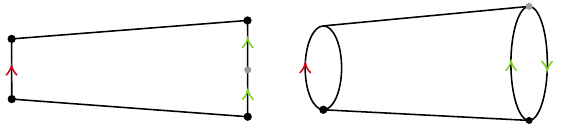}
    \caption{A relator and tube corresponding to some $t_j$.}
\end{figure}

\begin{itemize}
    \item If $j$ lies in a clean component of $\Gamma$, then $n_g/k = (2W+KLN^T+6C+1)|p_g|/k = (2W+KLN^T+6C+1)|p_{h}|/l = n_{h}/l$ is an integer by definition of the potential. Construct the tube $T_j$ from $n_g/k=n_h/l$ copies of the $2$-cell corresponding to the relator $t_jw_{g}^k t_j^{-1}w_{h}^{-l}$ glued along their $1$-cells labeled $t_j$ to form an annulus. This tube has two boundary cycles, labeled by $w_g^{n_g}$ and $w_{h}^{n_{h}}$. The word $w_g$ has $m=d/n_g$ distinct elevations to $V^{(1)}$ while $w_{h}$ has $n=d/n_{h}$ distinct elevations. Glue, by label-preserving homeomorphisms, the ends of $k$ copies of $T_j$ along each elevation of $w_g$, rotating the attaching map of each copy so that the endpoints of the $1$-cells labeled by $t_j$ are shifted by $|w_g|$ with each new copy. This provides a total of $mk = nl$ tubes, each having one boundary component glued to $V^{(1)}$. Glue the remaining boundary components along elevations of $w_h$ using label-preserving homeomorphisms, with $l$ tubes arriving at each elevation, and as before rotating the attaching maps by $|w_h|$ from one tube to the next (it only matters here that $l$ tubes arrive at each elevation, the matching of tubes and elevations can be done arbitrarily).

    \begin{figure}[h]
        \centering
        \includegraphics[width=0.5\textwidth]{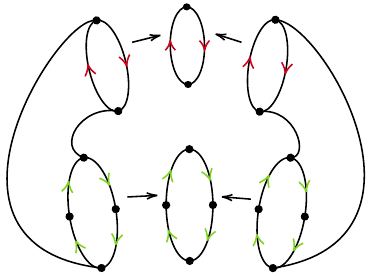}
        \caption{Tubes when $t_j$ belongs to a clean component. The $t_j$ edges with blue arrows are identified to form a tube. The attaching maps of the tubes are all homeomorphisms.}
    \end{figure}

    \item If $j$ lies in a non-clean component of $\Gamma$, one of $k,l$ is equal to $\pm 1$ by Proposition~\ref{nonclean2}, and $n_g=n_h$ by definition. Assume $|k|= 1$; the construction for $|l|=1$ is symmetric. Construct the tube $T_j$ from $n_g=n_h$ copies of the $2$-cell corresponding to the relator $t_jw_{g}^k t_j^{-1}w_{h}^{-l}$ glued along their $1$-cells labeled $t_j$ to form an annulus. This tube has two boundary cycles labeled by $w_g^{kn_g}$ and $w_h^{ln_h}$. The words $w_g$ and $w_h$ have the same number $m = d/n_g=d/n_h$ of elevations to $V^{(1)}$. Arbitrarily pair these elevations. For each pair, glue the boundary components of a copy of $T_j$ along the two elevations so that the gluing is a label-preserving degree $|l|$ cover on the $h$ side and a label-preserving homeomorphism on the $g$ side.

    \begin{figure}[h]
        \centering
        \includegraphics[width=0.3\textwidth]{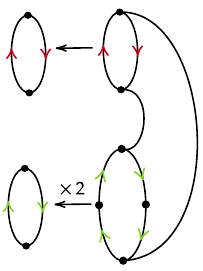}
        \caption{A tube when $t_j$ belongs to a non-clean component. One gluing is degree 2 and the other a homeomorphism.}
    \end{figure}
\end{itemize}

The space $X$ obtained by making the above construction for each edge $j$ in $\Gamma$ has a graph of spaces structure with a single vertex space $V^{(1)}$ and one edge space per tube. Let $\pi\colon X\to \GG$ be the projection to the underlying graph, a bouquet of circles. See Figure~\ref{fig:bigDiagram}.
\end{defi}

\begin{rem}
\label{rem:stableNotInterior}
Let $g\in A$ and  $e_g$ be some elevation of $w_g$ to $V^{(1)}$. Note that the decomposition of $e_g$ as a concatenation of subpaths labeled $w_g$ is unique, otherwise $w_g$ would be a proper power.

Assume that $g$ is the terminal endpoint of $j$ and $k$, two (not necessarily distinct) edges of $\Gamma$. Let $x_j$, $x_k$ be vertices of the image of $e_g$ that are terminal endpoints of $1$-cells of $X$ labeled by $t_j$ and $t_k$ respectively. By construction of $X$, each of these two $1$-cells belongs to a unique tube attached to $e_g$. Moreover, each of these two attaching maps at $e_g$ defines a decomposition of $e_g$ into subpaths labeled $w_g$, one starting at $x_j$ the other starting at $x_k$ respectively. By the above observation, the two decompositions coincide, and every subpath with endpoints $x_j$ and $x_k$ immersed in the image of $e_g$ must be labeled by some power of $w_g$.

The same holds symmetrically in the three following situations:
\begin{itemize}
    \item $g$ is the initial endpoint of $j$, and $x_j$ is the initial endpoint of a $1$-cell labeled $t_j$ instead
    \item $g$ is the initial endpoint of $k$ and $x_k$ is the initial endpoint of a $1$-cell labeled $t_k$ instead
    \item both of the above
\end{itemize}
\end{rem}

\begin{lm}
The cover $V^{(1)}\to (P_{\overline F})^{(1)}=P_F$ extends to a degree $d$ cover $X\to P_G$.
\end{lm}
\begin{proof}
$P_F$ is a subcomplex of $P_G$, and there is a unique $0$-cell $x$ in both complexes. Additionally, $V^{(1)}$ is a subcomplex of $X$ with the same $0$-skeleton. Since $X$ was constructed from $V^{(1)}$ by adding tubes, the additional $1$-cells of $X$ are oriented and labeled with stable letters $t_j$. Since $P_G$ has a single $0$-cell and oriented loops for each stable letter, there is a map $X^{(1)}\to P_G$ extending $V^{(1)}\to P_F\subset P_G$ respecting the orientation and labels of the additional $1$-cells. Since each $2$-cell of $X$ is attached by a cycle representing a relator of $G$, the map $X^{(1)}\to P_G$ extends to $X\to P_G$ by mapping each $2$-cell of $X$ to the unique $2$-cell in $P_G$ whose attaching map has the same label. 

$P_G$ has a $1$-cell for each $a_i$, $i\in I$ and $t_j$, $j\in J$. Thus there are $2(|I|+|J|)$ ends of $1$-cells incident at $x$, and $\lk(x)$ has vertices which can be denoted $a_i^\pm$, $i\in I$, and $t_j^\pm$, $j\in J$. For each $g\in A$, let $(w_g)_i$ denote the $i$-th letter of $w_g$, for $i=1,\dots, |w_g|$. If there exists a relator of the form $t_jw_g^kt_j^{-1}=w_{h}^l$, then in $\lk(x)$ the vertex $t_j^-$ is joined by edges to $(w_g)_1^-$ and $(w_g)_{|w_g|}^+$; and $t_j^+$ is joined by edges to $(w_h)_1^-$ and $(w_h)_{|w_{h}|}^+$. Additionally, for $i=1,\dots, |w_g|$, there are $k$ edges joining the two vertices $(w_g)_i^+$ and $(w_g)_{i+1}^-$. There are also $k-1$ edges joining $(w_g)_{|w_g|}^+$ and $(w_g)_1^-$ (these two vertices are also joined by a path of length $2$ through $t_j^-$). Likewise, for $i=1,\dots, |w_{h}|$, there are $l$ edges joining $(w_{h})_i^+$ to $(w_{h})_{i+1}^-$ and $l-1$ edges joining $(w_{h})_{|w_{h}|}^+$ and $(w_{h})_1^-$ (these two vertices are also joined by a path of length $2$ through $t_j^+$). All edges of $\lk(x)$ are of the above form for some relator. 

Let $\hat x\in X$ be a $0$-cell, and consider the map $\lk(\hat x)\to \lk(x)$ induced by $X\to P_G$. For $g\in A$, let $w_{g,i}$ be the cyclic permutation of $w_g$ beginning with the $i$-th letter for $i=1,\dots, |w_g|$. By Remark~\ref{rem:stableNotInterior}, elevations $e_{g,i}$ and $e_{g,j}$ of $w_{g,i}$ and $w_{g,j}$ beginning at $\hat x$ have the same image if and only if $i=j$. The elevation $e_{g,i}$ is an embedded cycle with label $w_{g,i}^{n_g}$ beginning at $\hat x$. Consider a relator $t_jw_g^kt_j^{-1}w_{h}^{-l}$ of $G$ where $w_g$ appears. We separately consider the cases of $j$ belonging to a clean component and a non-clean component of $\Gamma$.

\begin{itemize}
    \item If $j$ belongs to a clean component of $\Gamma$, then there are $k$ copies of the tube $T_j$ attached by homeomorphisms to each elevation $e_{g,i}$, whose attaching maps are rotated by $|w_g|$ from one copy to the next. See Definition~\ref{defi:tubeConstruction}. Of the copies of $T_j$ attached to the elevation $e_{g,1}$, exactly one has a $1$-cell labeled $t_j$ incident to $\hat x$ on its $t_j^-$ end, contributing a path of length $2$ joining $(w_g)_{|w_g|}^+$ to $t_j^-$ to $(w_g)_{1}^-$ in $\lk(\hat x)$. The remaining $(k-1)$ copies of $T_j$ attached to $e_{g,1}$ each contribute an edge joining $(w_g)_{|w_g|}^+$ to $(w_g)_1^-$ in $\lk(\hat x)$. Likewise, for each $i=2,\dots, |w_g|$, the $k$ copies of $T_j$ attached to $e_{g,i}$ contribute $k$ edges joining $(w_g)_{i-1}^+$ to $(w_g)_{i}^-$. A similar consideration of the elevations of $w_h$ shows that every edge of $\lk(x)$ coming from a corner of the $2$-cell $R$ of the relator $t_jw_g^kt_j^{-1}=w_{h}^l$ is the image of a unique edge of $\lk(\hat x)$ coming from a corner of a $2$-cell mapping to $R$.

    \item If $j$ belongs to a non-clean component of $\Gamma$, one of $k,l$ appearing in the relator $t_jw_g^kt_j^{-1}=w_{h}^l$ is equal to $\pm 1$ by Proposition~\ref{nonclean2}. Assume $|k|=1$; the case $|l| = 1$ is symmetric. For each $i=1,\dots, |w_g|$, there is a copy of the tube $T_j$ attached by a homeomorphism to $e_{g,i}$. The copy of $T_j$ attached to $e_{g,1}$ contributes a path of length $2$ joining $(w_g)_{|w_g|}^+$ to $t_j^+$ to $(w_g)_1^-$, and the other copies of $T_j$ contribute an edge joining $(w_g)_{i-1}^+$ to $(w_g)_i^-$ for each $i=2,\dots, |w_g|$. For each $i=1,\dots, |w_{h}|$ there is a copy of $T_j$ attached by a degree $l$ covering map to $e_{h,i}$. The $n_h$ $2$-cells of $T_j$ decompose one of its boundary cycles into a concatenation of $n_h$ subpaths labeled $w_h^l$, with a $t_j^-$ end of a $1$-cell labeled $t_j$ incident at each concatenation point. In the copy of $T_j$ glued to $e_{h,i}$, this boundary cycle is attached to the image of $e_{h,i}$ via a degree $l$ label-preserving covering map. Since $n_{h}$ is coprime to $l$, the $n_h$-many $t_j^-$ ends incident to $e_{h,i}$ coming from this copy of $T_j$ are evenly spaced $|w_h|$ apart in $e_{h,i}$. More precisely, each subpath of $e_{h,i}$ joining two closest incident $t_j^-$ ends is labeled $w_h$. In particular, there is a single $t_j^{-}$ end incident at $\hat x$, coming from the copy of $T_j$ glued to the image of $e_{h,1}$. Since this gluing is a degree $l$ covering map, this copy of $T_j$ contributes a path of length $2$ joining $(w_{h})_{|w_{h}|}^+$ to $t_j^-$ to $(w_{h})_1^-$ and $l-1$ edges joining $(w_{h})_{|w_{h}|}^+$ to $(w_{h})_1^-$ in $\lk(\hat x)$. Likewise, each copy of $T_j$ attached to some elevation $e_{h,i}$ for $i=2,\dots, |w_{h}|$ contributes $l$-many edges joining $(w_{h})_{i-1}^+$ to $(w_{h})_i^-$. Once again, every edge of $\lk(x)$ coming from a corner of the $2$-cell $R$ of the relator $t_jw_g^kt_j^{-1}=w_{h}^l$ is the image of a unique edge of $\lk(\hat x)$ coming from a corner of a $2$-cell mapping to $R$.
\end{itemize}

Since $V^{(1)}\to P_F$ is a covering map restricting $X\to P_G$, the map $\lk(\hat x)\to \lk(x)$ induces a bijection between vertices of $\lk(\hat x)$ coming from $1$-cells of $V^{(1)}$ and vertices of $\lk(x)$ coming from $1$-cells of $P_F$.

Since every edge of $\lk(x)$ comes from a corner of some $2$-cell $R$ of $P_G$, it is the image of an edge of $\lk(\hat x)$ by the above arguments. This proves surjectivity of $\lk(\hat x)\to \lk(x)$ on edges. By construction of $P_G$, every $1$-cell labeled by a stable letter is in the boundary of a $2$-cell. Hence, the corresponding vertices of $\lk(x)$ are not isolated. Therefore, every isolated vertex of $\lk(x)$ comes from a $1$-cell of $P_F$ and lies in the image of $\lk(\hat x)$. This proves that $\lk(\hat x)\to \lk(x)$ is surjective.

Furthermore, $1$-cells with distinct labels in $X$ are mapped to distinct $1$-cells in $P_G$. By the above arguments, each vertex $t_j^+$ or $t_j^-$ of $\lk(x)$ has a unique preimage in $\lk(\hat x)$. Hence $\lk(\hat x)\to \lk(x)$ is injective on vertices. Assume two edges $e$, $f$ of $\lk(\hat x)$ are mapped to the same edge of $\lk(x)$. The edges $e$ and $f$ correspond to corners of $2$-cells mapped to the same $2$-cell of $P_G$, which therefore must come from the same relator. By the uniqueness obtained above, $e=f$ and $\lk(\hat x)\to \lk(x)$ is injective on edges. By Lemma~\ref{lem:checkCover}, $X\to P_G$ is a degree $d$ covering map. 
\end{proof}

\begin{defi}
\label{def:theCover}
Consider $p$ a path in $X$ labeled by a representative of $w$. Such paths exist since they are lifts of paths in $P_G$ with the same label. The path $\pi \circ p$ in the underlying graph $\GG$ lifts to the universal cover $\widetilde{\GG}$ to a (possibly backtracking) path $\mathfrak p$. By residual finiteness of the free group $\pi_1\GG$, there exists an intermediate, regular, finite-degree cover $\widehat{\GG}\to \GG$ such that the map $\widetilde{\GG}\to \widehat{\GG}$ restricts to an embedding on the compact image of $\mathfrak{p}$. Let $\widehat X$ be the graph of spaces with underlying graph $\widehat{\GG}$ obtained from the graph of spaces $X$ and the cover $\widehat{\GG}\to \GG$ as in Remark~\ref{covergraph}. The vertex spaces of $\widehat X$ are isomorphic to $V^{(1)}$. Let $\hat \pi\colon \widehat X\to \widehat{\GG}$ be the projection to the underlying graph. See Figure~\ref{fig:bigDiagram}. By construction, the degree of the cover $\widehat X\to X$ is the same as the degree of $\widehat{\GG}\to \GG$. Thus the composition $\widehat X \to X\to P_G$ is a finite-degree covering.
\end{defi}

The cover $\widehat{X} \to P_G$ will satisfy the requirements of Proposition~\ref{mainprop}.

\begin{lm}
\label{returnfar}
Let $\hat q$ be a subpath of $\hat p$ whose label $s$ is an extremal subword of $w$. The following hold:
\begin{enumerate}
    \item The path $\hat \pi \circ \hat q$ fully backtracks. In particular, the endpoints $z_1,z_2$ of $\hat q$ belong to the same vertex space of $\widehat X$ (which like all such vertex spaces is isomorphic to $V^{(1)}$).
    \item There exists a unique $g\in A$ and a unique elevation $e_g$ of $w_g$ to this vertex space whose image contains $z_1$ and $z_2$. 
    \item Any path immersed in the image of $e_g$ joining $z_1$ and $z_2$ is longer than $K$. 
    \item A vertex excursion in $w$ adjacent to $s$ is not a power of $w_g$.
    \item For every subpath $\hat s$ of $\hat p$ containing $\hat q$, $\hat \pi\circ \hat s$ does not fully backtrack unless $\hat \pi\circ \hat s = \hat\pi\circ \hat q$.
\end{enumerate}

\end{lm}

\begin{proof}
Recall, by Corollary~\ref{cor:extremalnormalized}, we have assumed that every extremal subword of $w$ is normalized. In particular, $s$ decomposes as a concatenation $\tau w_h^l \tau^{-1}$ for some $h\in A$, $l\in \ZZ$, and $\tau$ a word in the $t_j^\pm$. The word $\tau$ corresponds to a path $\gamma$ in $\Gamma$ joining some vertex $g\in A$ to $h$. By Lemma~\ref{gammabacktrack} applied to $s$ and the backtrack $\gamma \gamma^{-1}$, the terminal end of $\gamma$ is not labeled $\pm 1$. The same lemma applied to $\tau$ and $\gamma$ shows every outgoing end of an edge of $\gamma$ is labeled $\pm 1$.

\emph{Claim:} Let $\tau'$ be a suffix of $\tau$, $\gamma'$ the corresponding subpath of $\gamma$, and $g'$ the starting point of $\gamma'$. Both endpoints $z_1'$, $z_2'$ of the subpath of $\hat q$ labeled $\tau'w_h^l\tau'^{-1}$ lie in the image of some elevation of $w_{g'}$ to a vertex space of $\widehat X$. Furthermore, if $|\tau'|>0$, the edge of $\hat q$ which starts at $z_1'$ and the edge which terminates at $z_2'$, which are labeled by some stable letter and its inverse, lie in the same tube (i.e.~edge space) of $\widehat X$.

Obtaining the claim for all suffixes of $\tau$, in particular for $\tau$ itself, proves Assertion~(1) and the existence part of Assertion~(2).

\begin{proof}[Proof of claim]
We work by induction over the length of $\tau'$. It is obvious when $\tau'$ and $\gamma'$ are trivial, with $g'=h$, since the subpath of $\hat q$ corresponding to $w_h^l$ lies in the image of some elevation by definition. Assume the claim holds for some suffix $\tau'$ and let $e_{g'}$ be the elevation of $w_{g'}$ containing $z_1'$ and $z_2'$. Let $\tau''=t_j\tau'$ be a larger suffix of $\tau$; the case where $\tau''=t_j^{-1}\tau'$ is symmetrical. Let $\gamma''$ be the corresponding subpath of $\gamma$, starting at $g''$. The first edge of $\gamma''$, labeled $t_j$, joins $g''$ to $g'$, has outgoing label $\pm 1$, incoming label $l\in \ZZ$, and lies in a non-clean component of $\Gamma$ as an edge of $\gamma$. By the structure of $X$ and its cover $\widehat X$, a single copy $\mathsf T$ of the tube $T_j$ is glued to the elevation $e_{g'}$ using a degree $|l|$ covering. Moreover, every edge labeled $t_j$ terminating in the image of $e_{g'}$ lies in $\mathsf T$. In particular, the edge of $\hat p$ joining $z_1''$ to $z_1'$ and the edge of $\hat p$ joining $z_2'$ to $z_2''$ lie in $\mathsf T$, proving the second part of the claim. Moreover, $z_1''$ and $z_2''$ both lie in the image of the elevation $e_{g''}$ of $w_{g''}$ on which one boundary component of $\mathsf T$ is glued, proving the first part of the claim.
\end{proof}

\bigskip

Now let $\hat r_g$ be a minimal length path immersed in the image of $e_g$ with endpoints $z_1$ and $z_2$. Since the first edge of $\gamma$ is incident to $g$ in $\Gamma$, by Remark~\ref{rem:stableNotInterior}, $\hat r_g$ is labeled by some power $w_g^k$. Hence, the word $sw_g^{-k} = \tau w_h^l \tau^{-1}w_g^{-k}$ labels the closed cycle $\hat q\hat r_g^{-1}$ in $\widehat X$, decomposing as two $\Gamma$-words $\tau w_h^l$, $\tau^{-1}w_g^{-k}$ corresponding to inverse paths $\gamma$, $\gamma^{-1}$ in $\Gamma$. By Lemma~\ref{normalize}, there exists a normalization of $sw_g^{-k}$ of the form $\tau w_h^m \tau^{-1}$ labeling a path $\hat q'$ homotopic rel. endpoints to $\hat q\hat r_g^{-1}$ in $\widehat X$. Therefore $\hat q'$ is closed and its subpath $\hat q''$ labeled by $w_h^m$ is closed and immersed in some vertex space. Letting $t_j^\epsilon$ be the first letter of $\tau$, the normalization process starts by replacing the subword $t_j^{-\epsilon} w_g^{-k}$ with some $w_{g'}^{k'}t_j^{-\epsilon}$, where $|k'|\leq N|k|$ by definition of $N$. Inductively, the subword $\tau^{-1} w_g^{-k}$ is replaced by $w_h^{l'}\tau^{-1}$, with $|l'|\leq N^{|\tau|}|k|\leq N^T|k|$, since $\tau$ appears as a subword of $w$, and $w$ contains $T$ stable letters. Finally, $|w_h^m| \leq |w_h^{l}|+|w_h^{l'}|\leq W + LN^T|k|$, since $w_h^l$ appears as a subword of $w$ and $|w_h|\leq L$.

Assuming that $|\hat r_g|\leq K$, then $|k|\leq K$ and $|w_h^m|\leq W+KLN^T$. Hence, $\hat q''$ corresponds to a closed path in $V^{(1)}$ with length at most $W+KLN^T$. By Remark~\ref{rem:shortLoopsLiftToLoops}, any lift $\widetilde{q''}$ to $\widetilde V=\widetilde{P_{\overline F}}$ is closed, immersed, and shorter than the boundary path of any $2$-cell, by definition of the $n_g$. Corollary~\ref{cor:shortLoopsAndPaths} provides a contradiction, proving $|\hat r_g|>K$ and Assertion~(3).

We now prove the uniqueness part of Assertion~(2). Assume $z_1$ and $z_2$ are contained in the image of some elevation $e_h$ of $w_h$ for some $h\in A$. Then there exist lifts $\tilde z_1$ and $\tilde z_2$ to $\widetilde V=\widetilde{P_{\overline F}}$ contained in lifts of $e_g$ and $e_h$, which are boundary cycles of $2$-cells, labeled by $w_g^{n_g}$ and $w_h^{n_g}$. Corollary~\ref{cor:relatorsEmbedded} implies that the intersection of these boundary cycles is convex. If the cycles are distinct, then Lemma~\ref{lem:boundedPiece} implies the intersection has diameter at most $C$. But we have shown that $z_1$ and $z_2$ are at distance at least $K>C$ in the image of $e_g$. Thus the cycles are not distinct, hence $g=h$ and $e_g=e_h$.

Since $g$ is the starting point of $\gamma$, a subword of $w$ of the form $w_g^n s$ or $s w_g^n$ with $n\neq 0$ would contradict the maximality of the extremal subword $s$, proving (4).

\begin{figure}[h]
    \centering
    \includegraphics[width=0.5\textwidth]{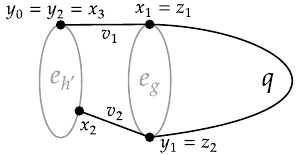}
    \caption{Using Lemma~\ref{lem:butterflyBalloon} to prove (5).}
    \label{fig:noExtension}
\end{figure}

Finally, let $\hat s$ be a subpath of $\hat p$ containing $\hat q$, and assume $\hat \pi\circ\hat s$ fully backtracks and properly contains $\hat \pi\circ \hat q$. It suffices to consider the minimal case where the label of $\hat s$ decomposes as a concatenation of the form $t_{j'}^{\pm 1}v_1\tau w_h^l\tau^{-1} v_2 t_{j'}^{\mp 1}$, with $v_1$ and $v_2$ two (possibly trivial) vertex excursions, and $j'\in J$. Let $y_2$ be the second vertex of $\hat s$ and $x_2$ the second-to-last vertex of $\hat s$. The first and last edge of $\hat s$ project to the same edge of $\GG$, hence lie in the same tube. Thus, the vertices $x_2$ and $y_2$ belong to some elevation $e_{h'}$ of some $w_{h'}$, $h'\in A$. Let $y_0 = y_2$, $x_1 = z_1$, $y_1 = z_2$, and $x_3 = y_2$. See Figure~\ref{fig:noExtension}. Apply Lemma~\ref{lem:butterflyBalloon} with $k=2$, $e_1 = e_{g}$ and $e_2 = e_{h'}$, $(x_1,x_2,x_3)$ and $(y_0,y_1,y_2)$. Assumption~3 holds because $\dist(y_0,x_1)+\dist(y_1,x_2)+\dist(y_2,x_3) \leq |v_1|+|v_2|+0\leq W$. Assumption~2 holds because $\dist(x_1,y_1)>K>W+3C$ and $\dist(x_2,y_2)>\dist(x_1,y_1)-(\dist(y_0,x_1)+\dist(y_1,x_2)) > K - W>W+3C$. Assumption~4a holds because of the choice $y_0 = x_3$. Therefore, $e_{g}=e_{h'}$. Let $\hat s_1$ (resp. $\hat s_2$) be the subpaths of $\hat s$ joining $y_0$ to $x_1$ (resp. $y_1$ to $x_2$), so that $\hat s_1$ is labeled $v_1$ and $\hat s_2$ is labeled $v_2$. By Lemma~\ref{lem:shortinelevation}, since $|v_1|,|v_2|\leq W$, both $\hat s_1$ and $\hat s_2$ are subpaths of $e_g = e_{h'}$. By definition of $\hat s$, $y_0$ and $x_1$ are the endpoints of edges labeled by stable letters ($t_{j'}$ and the first letter of $\tau$ respectively) incident to $g$ in $\Gamma$ with the correct orientations to apply Remark~\ref{rem:stableNotInterior}. By this remark, $\hat s_1$ is a power of $w_g$, and so is $\hat s_2$ by the same argument. Thus, the label of $\hat s$ cannot properly contain $s$ by definition of extremal subwords, a contradiction. Therefore, either $\hat \pi\circ \hat s = \hat \pi\circ \hat q$ or $\hat \pi\circ \hat s$ does not fully backtrack, proving (5).
\end{proof}

For the next lemma, recall a maximal backtrack is a fully backtracking (hence non-trivial) subpath that is maximal under inclusion given a fixed midpoint.
\begin{lm}
\label{lem:extremal}
Let $\sigma\sigma^{-1}$ be a maximal backtrack of $\hat \pi\circ \hat p$. There exists an extremal subword $s$ of $w$, corresponding to a subpath $\hat q$ of $\hat p$, such that  $\sigma\sigma^{-1} =\hat \pi\circ \hat q$.
\end{lm}

\begin{proof}
Let $f_1$, $f_2$ be the edges of $\hat p$ mapped by $\pi$ to the last edge of $\sigma$ and the first edge of $\sigma^{-1}$ respectively. Assume that $f_1$ is labeled by a stable letter $t_j$; the case where $f_1$ is labeled by the inverse $t_j^{-1}$ is symmetrical. Then $f_2$ is labeled by $t_j^{-1}$, and the subpath $\hat r$ of $\hat p$ between $f_1$ and $f_2$ is labeled by a vertex excursion $u$ of $w$. Let $t_jw_g^{k}t_j^{-1} = w_h^l$ be the relation in $G$ corresponding to $t_j$. Since $f_1$ and $f_2$ project to the same edge of $\widehat \GG$, they are in the same edge space of $\widehat X$, a copy of the tube $T_j$, and both endpoints of $\hat r$ belong to the same elevation $e_g$ of $w_g$ on which this tube is attached. By Lemma~\ref{lem:shortinelevation}, since $|\hat r|\leq W$, $\hat r$ is in fact a subpath of $e_g$. Moreover, both endpoints of $\hat r$ are endpoints of edges labeled $t_j$ oriented towards $e_g$. By Remark~\ref{rem:stableNotInterior}, the label $u$ of $\hat r$ is a power of $w_g$, i.e., the subword $t_jut_j^{-1}$ of $w$ is a $\Gamma$-word.

If $j$ belongs to a clean component of $\Gamma$, the gluing of the copy of $T_j$ at $e_g$ is by a homeomorphism. Then the subpath $f_1\hat rf_2$ of $\hat p$ is homotopic rel. endpoints inside this tube to a path remaining inside a boundary component of the tube. The label of the path obtained from $\hat p$ by applying this homotopy to the subpath $f_1\hat rf_2$ represents the same element of $G$ as $w$, and this contradicts the minimality of the number of stable letters of $w$. Therefore, $j$ belongs to a non-clean component of $\Gamma$, and $t_jut_j^{-1}$ is a non-clean $\Gamma$-word. This subword of $w$ is contained in some extremal subword with apex $u$. Since all extremal subwords of $w$ are normalized, this extremal subword is of the form, $\tau u \tau^{-1}$, where $\tau$ is a word in the $t_{j}^\pm$. Let $\hat q$ be the subpath of $\hat p$ corresponding to this extremal subword. By Lemma~\ref{returnfar}~(1), $\hat \pi \circ \hat q$ fully backtracks with the same midpoint $\hat\pi (u)\in \widehat\GG$ as $\sigma\sigma^{-1}$. Since $\sigma \sigma^{-1}$ is a maximal backtrack, $\hat \pi\circ \hat q$ is a subpath of $\sigma \sigma^{-1}$. By Lemma~\ref{returnfar}~(5), this subpath cannot be proper.
\end{proof}

\begin{lm}
\label{lem:gammaatmostonebacktrack}
    Let $\hat q$ be a subpath of $\hat p$ labeled by a non-clean $\Gamma$-word. Then $\hat \pi \circ \hat q$ backtracks at most once. Consequently, two distinct maximal backtracks of $\hat \pi\circ \hat p$ have disjoint interiors in the domain.
\end{lm}
\begin{proof}
    Let $s$ be the label of $\hat q$, and let $\gamma$ be the corresponding path in $\Gamma$. The midpoint of a backtrack of $\hat \pi \circ \hat q$ corresponds to a vertex excursion $v$ preceded and followed in $s$ by inverse stable letters. The excursion $v$ is a power of some $w_g$, for $g$ a vertex in some non-clean component of $\Gamma$. The path $\gamma$ in $\Gamma$ then backtracks at $g$. Assuming $\hat \pi \circ \hat q$ contains two backtracks, $\gamma$ contains two backtracks as well, at vertices $g_1$ and $g_2$. Let $\gamma'$ be the non-trivial subpath of $\gamma$ delimited by the midpoints $g_1$ and $g_2$ of those two backtracks. By Lemma~\ref{gammabacktrack} applied with $s$ and $\gamma$ neither end of $\gamma'$ is labeled $\pm 1$. Since $\gamma'$ lies in a non-clean component, by Lemma~\ref{gammabacktrack} again this is a contradiction. Therefore $\hat \pi \circ \hat q$ backtracks at most once.

    Let $\sigma_1\sigma_1^{-1}$, $\sigma_2\sigma_2^{-1}$ be two distinct maximal backtracks of $\hat \pi\circ \hat p$. By Lemma~\ref{lem:extremal}, there exist extremal subwords $s_1$, $s_2$ of $w$ labeling subpaths of $\hat p$ projecting via $\hat \pi$ to $\sigma_1\sigma_1^{-1}$ and $\sigma_2\sigma_2^{-1}$. If the interiors of $\sigma_1\sigma_1^{-1}$ and $\sigma_2\sigma_2^{-1}$ overlap in the domain, then $s_1$ and $s_2$ overlap on at least one stable letter. Therefore, their union $s$ is a non-clean $\Gamma$-word. Let $\hat q$ be the subpath of $\hat p$ labeled by $s$. Then $\hat \pi\circ \hat q$ contains both $\sigma_1\sigma_1^{-1}$ and $\sigma_2\sigma_2^{-1}$, two distinct backtracks, a contradiction with the previous result. Thus, $\sigma_1\sigma_1^{-1}$ and $\sigma_2\sigma_2^{-1}$ have disjoint interiors
\end{proof}

\begin{lm}
\label{lem:onlyoneextremal}
    There does not exist a subpath $\hat s_1 \hat v\hat s_2$ of $\hat p$ where $\hat s_1$ and $\hat s_2$ are labeled by extremal subwords of $w$ and $\hat v$ is labeled by a vertex excursion of $w$, such that $\hat s_1$ and $\hat s_2$ have all their endpoints in the image of the same elevation $e_g$ of some $w_g$, $g\in A$. 
\end{lm}
\begin{proof}
    Let $s_1$, $v$, $s_2$ be the labels of $\hat s_1$, $\hat v$, $\hat s_2$ respectively, and assume the endpoints of these three paths lie in the image of $e_g$.
    By Lemma~\ref{lem:shortinelevation}, since $|v|\leq W$, $\hat v$ is a subpath of $e_g$. By Lemma~\ref{returnfar}~(1) and (2), the last edge of $\hat s_1$ and the first edge of $\hat s_2$ are labeled by stable letters whose corresponding edges in $\Gamma$ are incident at $g$, with the correct orientation to apply Remark~\ref{rem:stableNotInterior}. By this remark, $v$ is a power of $w_g$, hence $s_1vs_2$ is a $\Gamma$-word. It is non-clean since $s_1$ and $s_2$ are non-clean. Yet, the path $\hat \pi \circ (\hat s_1\hat v \hat s_2)$ backtracks twice at the apexes of $s_1$ and $s_2$, contradicting Lemma~\ref{lem:gammaatmostonebacktrack}.
\end{proof}

To conclude, we need a last combinatorial lemma:

\begin{defi}
\label{def:star}
Let $\mathfrak{p}$ be a path in a graph whose maximal backtracks have pairwise disjoint interiors in the domain. A \emph{star subpath} of $\mathfrak{p}$ is either a maximal backtrack or a maximal subpath decomposing as a concatenation of at least two backtracks.

Note that since every backtrack is contained in a maximal backtrack and maximal backtracks have disjoint interiors, in the latter case, all the backtracks involved in the concatenation are maximal and the decomposition is unique.
\end{defi}

\begin{lm}
\label{lem:tinsel}
Let $\mathfrak{p}$ be a non-trivial path in a tree $\mathcal{T}$. Assume the following:
\begin{itemize}
    \item $\mathfrak{p}$ does not decompose as a concatenation of full backtracks (i.e.~is not a star subpath of itself).
    \item Maximal backtracks in $\mathfrak{p}$ have disjoint interiors in the domain.
    \item For every pair of (unoriented) edges $e_1$, $e_2$ that immediately precede and follow a star subpath in $\mathfrak{p}$, $e_1\neq e_2$
\end{itemize}
Then $\mathfrak{p}$ is not closed.
\end{lm}
Note that the last assumption automatically holds when the star subpath consists of a single maximal backtrack.

\begin{proof}
First note that by definition of star subpaths and since maximal backtracks have disjoint interiors in the domain, distinct star subpaths of $\mathfrak{p}$ are entirely disjoint in the domain. Since $\mathfrak{p}$ is not a star subpath of itself, one of its edges is not contained in any star subpath.

Now consider $\overline{\mathfrak{p}}$ the non-trivial path obtained from $\mathfrak{p}$ by deleting all of its star subpaths (which are closed). The paths $\overline{\mathfrak{p}}$ and $\mathfrak{p}$ are homotopic rel. endpoints. In particular, they have the same endpoints.

Assume for sake of contradiction that $\overline{\mathfrak{p}}$ contains a backtrack $\mathfrak{q}\mathfrak{q}^{-1}$ with $|\mathfrak{q}|=1$. Then $\mathfrak{p}$ contains a subpath either of the form $\mathfrak{q}\mathfrak{q}^{-1}$ or $\mathfrak{q}\mathfrak{s}\mathfrak{q}^{-1}$, where $\mathfrak{s}$ is a star subpath. The first case is impossible since any backtrack $\mathfrak{q}\mathfrak{q}^{-1}$ in $\mathfrak{p}$ is contained in a star subpath, and all star subpaths were deleted to define $\overline{\mathfrak{p}}$. The second case is impossible because the edges preceding and following $\mathfrak{s}$ would be identical, contradicting a hypothesis. Therefore, $\overline{\mathfrak{p}}$ contains no backtracks, hence is immersed. Since $\overline{\mathfrak{p}}$ is a nontrivial, immersed path in a tree, its endpoints are distinct. Thus $\mathfrak{p}$ is not closed. 
\end{proof}

We can now finally move on to the proof of Proposition~\ref{mainprop}.

\begin{proof}[Proof of Proposition~\ref{mainprop}]
Recall $\hat p$ is a path labeled by $w$ in $\widehat X$. Assume that $\hat p$ is closed. Several cases arise, each yielding a contradiction.
\begin{itemize}
    \item Assume $\hat \pi\circ \hat p$ is constant. Then $\hat p$ is closed, immersed in a vertex space and of length $W$. By Remark~\ref{rem:shortLoopsLiftToLoops}, $w$ labels a closed path immersed in $\widetilde{P_{\overline F}}$. This path is shorter than any relator of $\overline F$, and Corollary~\ref{cor:shortLoopsAndPaths} yields a contradiction. 
    \item If $\hat \pi\circ \hat p=\sigma\sigma^{-1}$ is non-constant but fully backtracking, by Lemma~\ref{lem:extremal}, $w$ decomposes as a product $v_1sv_2$ where $v_1$, $v_2$ are vertex excursions and $s$ is extremal. By Lemma~\ref{returnfar}~(2), the endpoints of the subpath corresponding to $s$ lie in an elevation $e_g$ of $w_g$ for some $g\in A$. Since $w$ is closed, there is a path $\hat q$ labeled $v_1^{-1}v_2$ joining the endpoints of the subpath corresponding to $s$. By Lemma~\ref{lem:shortinelevation}, the immersed path homotopic in $V^{(1)}$ rel endpoints to $\hat q$ is embedded in the image of $e_g$, being no longer than $|\hat q| = |v_1|+|v_2|<W$. Since $|\hat q|<W<K$, this contradicts Lemma~\ref{returnfar}~(3). 

    \item Assume $\hat \pi\circ \hat p = \sigma_1\sigma_1^{-1}\cdots\sigma_k\sigma_k^{-1}$ decomposes as a concatenation of at least two backtracks. By Lemma~\ref{lem:gammaatmostonebacktrack}, maximal backtracks have disjoint interiors in $\hat\pi\circ \hat p$. Thus, as in Definition~\ref{def:star}, all the backtracks in the concatenation are maximal. By Lemma~\ref{lem:extremal}, $w$ decomposes as a concatenation $v_1s_1\cdots v_ks_kv_{k+1}$ where the $v_l$ are vertex excursions and the $s_l$ are extremal. Let $\hat v_1\hat s_1\cdots \hat v_k\hat s_k\hat v_{k+1}$ be the corresponding decomposition of $\hat p$ in subpaths. Note that $k\leq W$ since each $s_i$ is non-trivial. For $1\leq i\leq k$, let $x_i$, $y_i$ be the endpoints of $\hat s_i$. Let $y_0=y_k$, $x_{k+1}=y_k$, and apply Lemma~\ref{lem:butterflyBalloon}, the elevations $e_i$ being uniquely determined by Lemma~\ref{returnfar}~(2). Assumptions~1 and 2 hold by Lemma~\ref{returnfar}~(2) and (3), as $K>W+3C$. Assumption~3 holds because $\displaystyle \sum_{i=0}^k \dist(y_i,x_{i+1}) \leq \sum_{i=1}^{k+1} |v_{i}| \leq W$. Finally, Assumption~4a holds because $y_0 = x_{k+1}$ by definition. Thus there exists $i$ such that the endpoints of $\hat s_i$, $\hat v_{i+1}$, and $\hat s_{i+1}$ all lie in the image of the same elevation $e_i = e_{i+1}$. This contradicts Lemma~\ref{lem:onlyoneextremal}.

    \begin{figure}[h]
        \centering
        \includegraphics[width=0.6\textwidth]{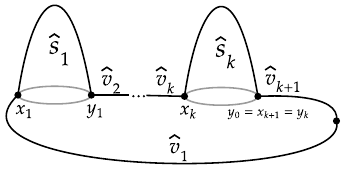}
        \caption{Using Lemma~\ref{lem:butterflyBalloon} to obtain a contradiction when $\hat\pi\circ\hat p$ is a concatenation of at least two full backtracks.}
    \end{figure}
    
    \item Otherwise, $\hat \pi\circ \hat p$ does not decompose as a concatenation of any number of backtracks. Hence, a lift $\mathfrak{p}$ of $\hat \pi\circ \hat p$ to the tree $\widetilde{\GG}$ is not a concatenation of backtracks either. By Lemma~\ref{lem:gammaatmostonebacktrack}, maximal backtracks in $\mathfrak{p}$ have disjoint interiors. Our goal is to apply Lemma~\ref{lem:tinsel} to $\mathfrak{p}$. We need to check its last assumption.
    
    Let $\mathfrak{p}'=\sigma_1\sigma_1^{-1}\cdots\sigma_k\sigma_k^{-1}$ be a star subpath of $\mathfrak{p}$, with $k\geq 2$ and assume the same edge precedes and follows $\mathfrak{p}'$ in $\mathfrak{p}$. The same edge precedes and follows the image of $\mathfrak{p}'$ in $\widehat \GG$ as well. Let $\hat p'$ be the longest subpath of $\hat p$ such that $\hat \pi\circ \hat p'$ coincides with the image of $\mathfrak{p}'$ in $\widehat \GG$. Then, the edges $f_1, f_2$ preceding and following $\hat p'$ in $\hat p$ project via $\hat \pi$ to the same edge. Hence $f_1$ and $f_2$ are labeled by the same stable letter, with opposite orientations, and lie in the same edge space of $\widehat X$. Letting $y_0$ and $x_{k+1}$ be the endpoints of $\hat p'$, this means that $y_0$ and $x_{k+1}$ belong to the same elevation of some $w_g$, $g\in A$. Moreover, by Lemma~\ref{lem:extremal}, $\hat p'$ decomposes as a concatenation $\hat v_1\hat s_1\cdots \hat v_{k}\hat s_{k}\hat v_{k+1}$ where each $\hat s_i$ is labeled by an extremal subword of $w$, and each $\hat v_i$ is labeled by a vertex excursion. The $\hat s_i$ are non-trivial, hence $k\leq W$. For $1\leq i\leq k$, let $x_i, y_i$ be the endpoints of $\hat s_i$. Note that the endpoints of $\hat v_i$ are $y_{i-1}, x_i$ for every $1\leq i\leq {k+1}$. Apply Lemma~\ref{lem:butterflyBalloon} this collection of points, the elevations $e_i$ being uniquely determined by Lemma~\ref{returnfar}~(2). Assumptions~1, 2, and 3 hold exactly as in the previous case, and Assumption~4b holds as well. Once again, there exists $i$ such that the endpoints of $\hat s_i$, $\hat v_{i+1}$, and $\hat s_{i+1}$ all lie in the image of the same elevation $e_i = e_{i+1}$. This is a contradiction by Lemma~\ref{lem:onlyoneextremal}, hence $\mathfrak{p}'$ cannot be preceded and followed by the same edge. By Lemma~\ref{lem:tinsel}, $\mathfrak{p}$ is not closed in the tree $\widetilde{\GG}$. However, the covering $\widetilde{\GG}\to \widehat{\GG}$ was chosen to be injective on the image of $\mathfrak{p}$. Thus the projection $\hat \pi\circ \hat p$ of $\mathfrak{p}$ to $\widehat \GG$ is non-closed as well. Therefore, $\hat p$ cannot be closed, the final contradiction.

    \begin{figure}[h]
        \centering
        \includegraphics[width=0.6\textwidth]{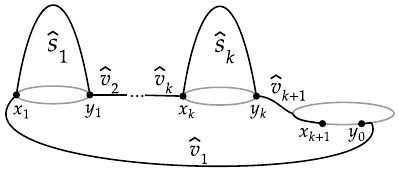}
        \caption{Using Lemma~\ref{lem:butterflyBalloon} to obtain a contradiction when $\hat\pi\circ\hat p$ is nontrivial and not a concatenation of full backtracks.}
    \end{figure}
\end{itemize}
\end{proof}

It is a clear consequence of Proposition~\ref{mainprop} that $G$ is residually finite. Recall that $G$ was taken to split as a finite multiple HNN-extension of a finite-rank free group, with cyclic edge groups. However, by Remark~\ref{rem:canUseHNN}, the same conclusion will hold for any finite graph of (arbitrary rank) free groups with cyclic edge groups. Together with Proposition~\ref{vuimpliesbs} and Remark~\ref{rem:algorithm}, we have proved:

\begin{thm}[Theorem~\ref{mthm:A}]
A group splitting as a finite graph of free groups with cyclic edge groups is residually finite if and only if all of its Baumslag-Solitar subgroups are residually finite. Morevoer, the residual finiteness can be determined algorithmically given a presentation of the splitting.
\end{thm}

\begin{rem}
The "only if" part of the theorem, given by Proposition~\ref{vuimpliesbs}, still holds without the finiteness assumption on the splitting, but the "if" direction fails. Indeed, consider the following standard presentation for a multiple HNN-extension.
\[G = \gen{a,b, (t_n)_{n\geq 1}\mid  t_n a t_n^{-1} = (ab^n)^{n+1}}\]

The element $a$ has roots of all orders, thus lies in the kernel of every map from $G$ to a finite group. Hence $G$ is not residually finite. The graph $\Gamma$ for the splitting is a depth $1$ infinite rooted tree with $a$ as the root and the $ab^n$, $n\geq 1$ as the leaves. Note that conjugates of edge groups pairwise intersect trivially in $\gen{a,b}$ except when they are both the same conjuate of the edge group $\gen{a}$.

Assume $g\in G$ is very unbalanced for $h,m,n$. By Lemma~\ref{unbalancedelliptic}, $g$ acts elliptically and $h$ acts loxodromically on the Bass-Serre tree of the HNN-extension. Let $v$ be a vertex fixed by $g$. Then $g^n = hg^mh^{-1}$ fixes the geodesic joining $v$ and $hv\neq v$. Since conjugates of edge groups are pairwise disjoint except when they are both the same conjugate of $\gen{a}$, the vertices $v$ and $hv$ are joined by a path of length at most $2$ stabilized by $g$ and there exists $h'\in G$ and $p,q$ non-negative distinct integers satisfying the following: $h'g^mh'^{-1}$ is a power of $ab^p$ and $h'g^nh'^{-1}$ is a power of $ab^q$. Thus $ab^p$ and $ab^q$ share a common power, a contradiction since $p\neq q$ and $\gen{a,b}$ embeds in $G$. Therefore, $G$ has no very unbalanced elements.\end{rem}

Our work leaves the following questions unanswered.

\begin{q}
Let $G$ be a group splitting as an infinite graph of free groups with cyclic edge groups. Assume $G$ has no very unbalanced elements nor non-trivial elements with roots of all orders. Is $G$ residually finite?
\end{q}

\begin{q}
Is there a combinatorial property of the graph $\Gamma$ encoding when $G$ is linear? residually $p$-finite for a prime $p$?
\end{q}


We also expect an argument very similar to ours to characterize residual finiteness in graphs of hyperbolic special groups with cyclic edge groups.
\begin{conj}
A group splitting as a finite graph of special groups with cyclic edge groups is residually finite if and only if it contains no very unbalanced elements.
\end{conj}
The proof of this conjecture would rely on the Malnormal Special Quotient Theorem, as well as cubical small-cancellation.
\bibliographystyle{alpha}
{\footnotesize
\bibliography{bibliography}}

\end{document}